# NEEDLES AND STRAW IN HAYSTACKS: EMPIRICAL BAYES ESTIMATES OF POSSIBLY SPARSE SEQUENCES

By Iain M. Johnstone[1] and Bernard W. Silverman[2]

*Stanford University and Oxford University*

An empirical Bayes approach to the estimation of possibly sparse sequences observed in Gaussian white noise is set out and investigated. The prior considered is a mixture of an atom of probability at zero and a heavy-tailed density $\gamma$, with the mixing weight chosen by marginal maximum likelihood, in the hope of adapting between sparse and dense sequences. If estimation is then carried out using the posterior median, this is a random thresholding procedure. Other thresholding rules employing the same threshold can also be used. Probability bounds on the threshold chosen by the marginal maximum likelihood approach lead to overall risk bounds over classes of signal sequences of length $n$, allowing for sparsity of various kinds and degrees. The signal classes considered are "nearly black" sequences where only a proportion $\eta$ is allowed to be nonzero, and sequences with normalized $\ell_p$ norm bounded by $\eta$, for $\eta > 0$ and $0 < p \leq 2$. Estimation error is measured by mean $q$th power loss, for $0 < q \leq 2$. For all the classes considered, and for all $q$ in $(0, 2]$, the method achieves the optimal estimation rate as $n \to \infty$ and $\eta \to 0$ at various rates, and in this sense adapts automatically to the sparseness or otherwise of the underlying signal. In addition the risk is uniformly bounded over all signals. If the posterior mean is used as the estimator, the results still hold for $q > 1$. Simulations show excellent performance. For appropriately chosen functions $\gamma$, the method is computationally tractable and software is available. The extension to a modified thresholding method relevant to the estimation of very sparse sequences is also considered.

## 1. Introduction.

Received August 2002; revised August 2003.

[1]Supported by NSF Grant DMS-00-72661 and NIH Grants CA 72028 and R01 EB001988-08.

[2]Supported by the Engineering and Physical Sciences Research Council.

1.1. *Thresholding to find needles and straw.* There are many statistical problems where the object of interest is a high-dimensional parameter on which we have a single observation, perhaps after averaging, and subject to noise. Specifically, suppose that $X = (X_1, \ldots, X_n)$ are observations satisfying

(1) $$X_i = \mu_i + \epsilon_i,$$

where the $\epsilon_i$ are $N(0, 1)$ random variables, not too highly correlated. Let $\mu$ be the vector of means $\mu = (\mu_1, \mu_2, \ldots, \mu_n)$. Clearly, without some knowledge of the $\mu_i$ we are not going to be able to estimate them very effectively, and in this paper we consider the advantage that may be taken of possible sparsity in the sequence.

In what contexts do problems of this kind arise? Some examples are the following:

- In astronomical and other image processing contexts, the $X_i$ may be noisy observations of the pixels of an image, where it is known that a large number of the pixels may be zero.
- In the model selection context, there may be many different models that conceivably contribute to the observed data, but it is of interest to select a subset of the possible models. In this case, the individual $X_i$ are the raw estimates of the coefficients of the various models, renormalized to have variance 1.
- In data mining, we may observe many different aspects of an individual or population, and we are only interested in the possibly small number that are "really there"; this is much the same as the model selection situation, but couched in different language.
- In nonparametric function estimation using wavelets, the true wavelet coefficients at each level form a possibly sparse sequence, and the discrete wavelet transform yields a sequence of raw coefficients, which are observations of these coefficients subject to error. Wavelet approaches in nonparametric regression take advantage of this structure in a natural way. This context originally motivated the work of this paper but the potential applicability of the ideas developed is much wider.

A natural approach to all these problems is *thresholding*: if the absolute value of a particular $X_i$ exceeds some threshold $t$, then it is taken to correspond to a nonzero $\mu_i$ which is then estimated, most simply by $X_i$ itself. If $|X_i| < t$, then the coefficient $|\mu_i|$ is estimated to be zero. But how is $t$ to be chosen? The importance of choosing $t$ appropriately is illustrated by a simple example. Consider a sequence of 10,000 $\mu_i$, of which $m$ are nonzero and $(10{,}000 - m)$ zero. The nonzero values are allocated at random and are each generated from a uniform distribution on $(-5, 5)$. By varying the number $m$, sequences of different sparsities can be generated, as shown in Figure 1.

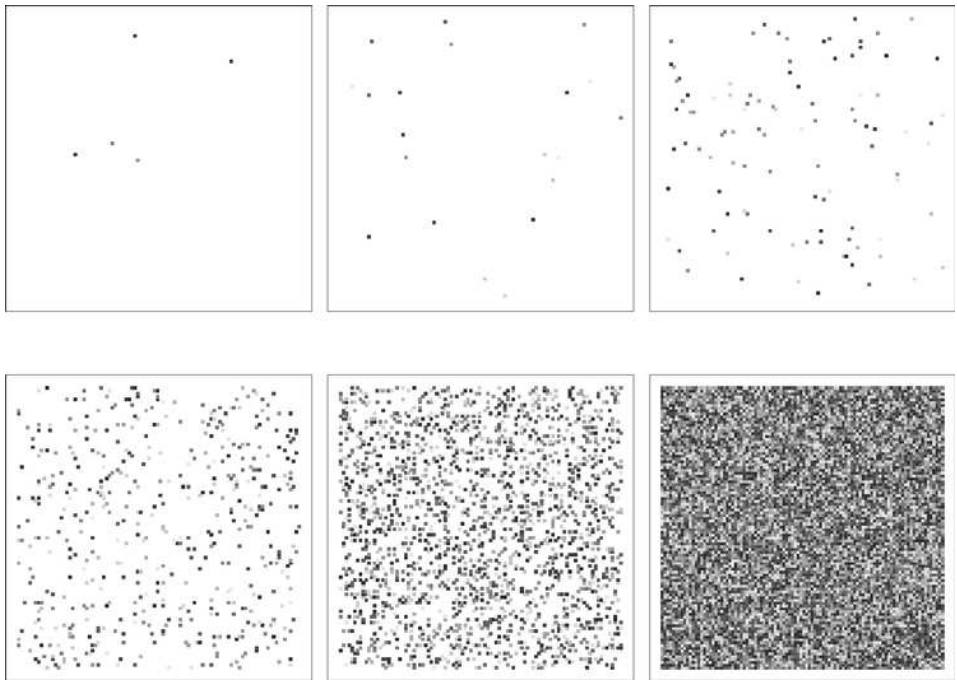

Fig. 1. *Absolute value of parameter images of various sparsity. Out of* 10,000 *pixels, the number of nonzero parameters is, from left to right:* 5, 20, 100 *in the top row and* 500, 2000, 10,000 *in the bottom row. Each nonzero parameter is chosen independently from a uniform distribution on* (−5, 5).

In this figure the 10,000 $\mu_i$ are arranged in a $100 \times 100$ pixel image. The absolute value of the image is plotted in gray scale in order to allow white to correspond to the value zero. Estimating a sparse signal is like finding needles in a haystack; it will be necessary to find which are the very few signal values that are nonzero, as well as to estimate them. On the other hand, estimating a dense signal is more like finding straw in a haystack; no longer will we be surprised if a particular $\mu_i$ is nonzero.

Independent Gaussian noise of variance 1 is added to the $\mu_i$ to yield a sequence $X_i$. The resulting images are shown in Figure 2. The average square estimation error yielded by thresholding $X_i$ with varying thresholds is plotted in Figure 3. Ignore the points marked by arrows for the moment. The number in the top right of each panel is the value of $m$, so $m = 5$ corresponds to a very sparse model, while $m = 10,000$ corresponds to a very dense model, with no zero parameter values at all. The naive estimator, estimating each $\mu_i$ by the corresponding $X_i$ without performing any thresholding at all, will produce an expected mean square error of 1. The scales in each panel are the same, and the threshold range is from 0 to $\sqrt{2 \log 10,000} \doteq 4.292$, the so-called *universal threshold* for a sample of this size.

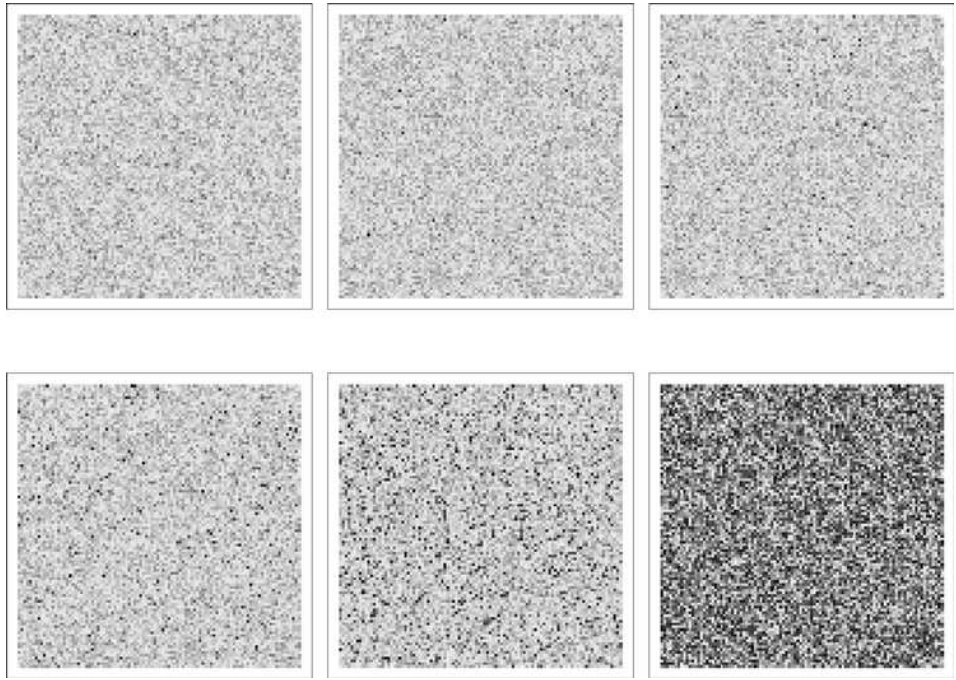

FIG. 2. *Absolute values of data $X_i$, result of adding Gaussian white noise to the images depicted in Figure* 1.

Three things can be seen from this figure. First, the potential gain from thresholding is very large if the true parameter space is sparse. For the sparsest signals considered in Figures 1 and 3, the minimum average square error achieved by a thresholding estimate is 0.01 or even less; see Figure 4 for a graph of minimum average square error against sparsity. Second, the appropriate threshold increases as the signal becomes more sparse. For the fully dense signal, no thresholding at all is appropriate, while for the sparsest signals, the best results are obtained using the universal threshold. Finally, it is important for the threshold to be tuned to the sparsity of the signal; if a threshold appropriate for dense signals is used on a sparse signal, or vice versa, the results are disastrous.

Thus, thresholding is a very promising approach, but the crucial aspect is the choice of threshold. A good threshold choice method will have several properties, as follows:

- It will be adaptive between sparse and dense signals, between finding "needles" and finding "straw."
- It will be stable to small changes in the data.
- It will be tractable to compute, with software available.
- It will perform well on simulated data and on real data.

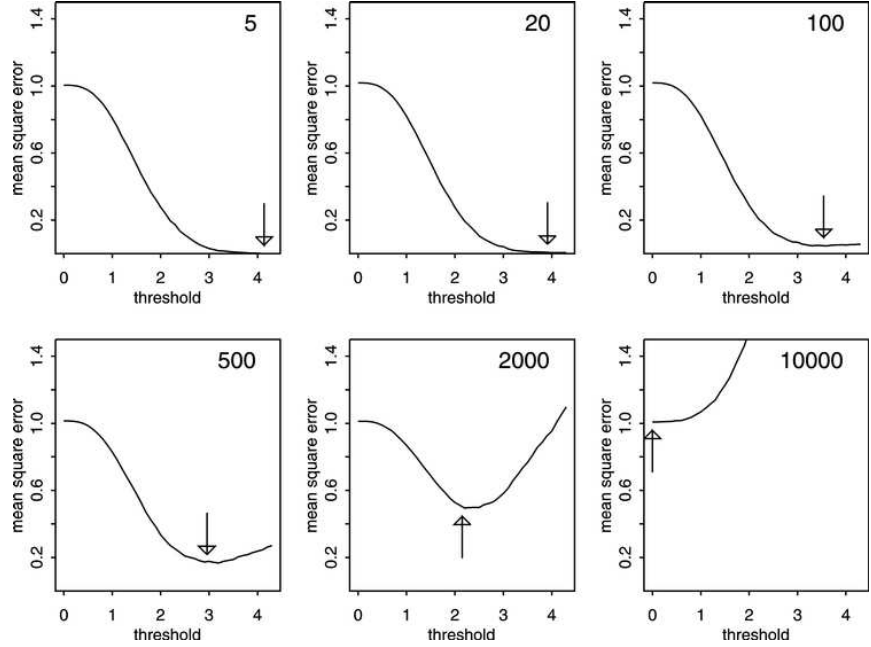

Fig. 3. *Mean square error of thresholding data obtained from the images in Figure 1 by adding Gaussian white noise. In each panel the arrow indicates the threshold chosen by the empirical Bayes approach. The prior used for the nonzero part of the distribution was a Laplace distribution with scale parameter $a = \frac{1}{2}$. Each plot is labeled by the number of nonzero pixels, out of 10,000, in the underlying signal.*

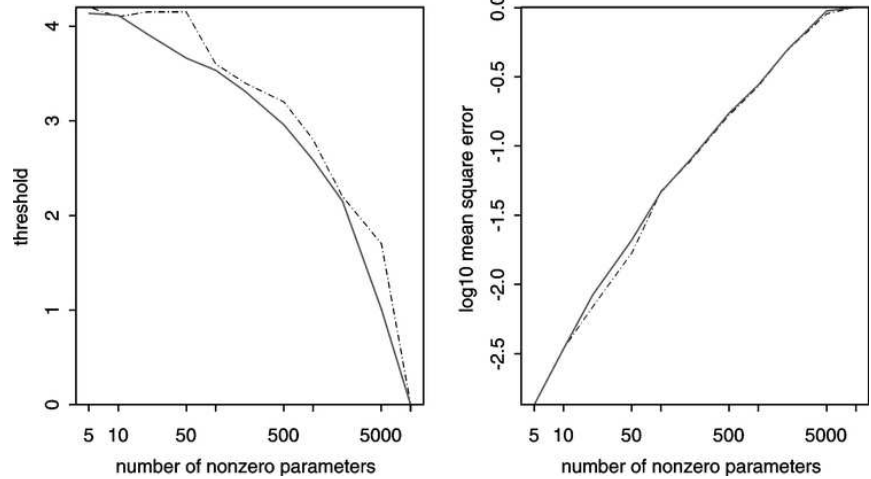

Fig. 4. *Left panel: threshold plotted against sparsity. The solid line is the threshold chosen by the empirical Bayes method, while the dashed line is the threshold that yields the minimum possible average square error. Right panel: log base 10 of the average square error yielded by the empirical Bayes threshold (solid line) and by the best possible threshold (dashed line). The models illustrated in Figure 1, and intermediate models, were used to construct these graphs.*

• It will have good theoretical properties.

In this paper we set out and investigate a fully automatic empirical Bayes thresholding method, which satisfies all these desiderata. In the example the method chooses the threshold values shown by the arrows in Figure 3. It can be seen that the empirical Bayes method is very good at tracking the minimum of the average square error. More details are given in Figure 4. The empirical Bayes thresholds are always close to the optimal thresholds, and—right across the range of sparsity considered—the average square error obtained by the empirical Bayes threshold is very close indeed to the best attainable average square error. A documented implementation *EbayesThresh* of our methodology in R and S-PLUS is available. See Johnstone and Silverman (2003) for details.

1.2. *Specifying the empirical Bayes method.*   In the present paper we concentrate attention on the case where the errors $\epsilon_i$ are independent. In some contexts this assumption is restrictive. While beyond the scope of the present paper, it is of obvious interest to extend our method and the supporting theory to dependent data, and this is a natural topic for future work.

The notion that many or most of the $\mu_i$ are near zero is captured by assuming that the elements $\mu_i$ have independent prior distributions each given by the mixture

$$(2) \qquad f_{\mathrm{prior}}(\mu) = (1 - w)\delta_0(\mu) + w\gamma(\mu).$$

The nonzero part of the prior, $\gamma$, is assumed to be a fixed unimodal symmetric density. In most previous work in the wavelet context mentioned above, the density $\gamma$ is a normal density, but we shall see that there are advantages in using a heavier-tailed prior, for example, a double exponential distribution or a distribution with tails that decay at polynomial rate.

For any particular value of the weight $w$, consider the posterior distribution of $\mu$ given $X = x$ under the assumption that $X \sim N(\mu, 1)$. Let $\hat{\mu}(x; w)$ be *median* of this distribution. For fixed $w < 1$, the function $\hat{\mu}(x; w)$ will be a monotonic function of $x$ with the thresholding property, in that there exists $t(w) > 0$ such that $\hat{\mu}(x; w) = 0$ if and only if $|x| \le t(w)$. Figure 5 shows the prior distribution and the posterior median function $\hat{\mu}(x; w)$ for the Laplace mixture prior with $a = 0.5$ and two different values of the weight $w$.

Let $g$ denote the convolution of the density $\gamma$ with the standard normal density $\phi$. The marginal density of the observations $X_i$ will then be

$$(1 - w)\phi(x) + wg(x).$$

We define the marginal maximum likelihood estimator $\hat{w}$ of $w$ to be the maximizer of the marginal log likelihood

$$\ell(w) = \sum_{i=1}^{n} \log\{(1 - w)\phi(X_i) + wg(X_i)\}$$

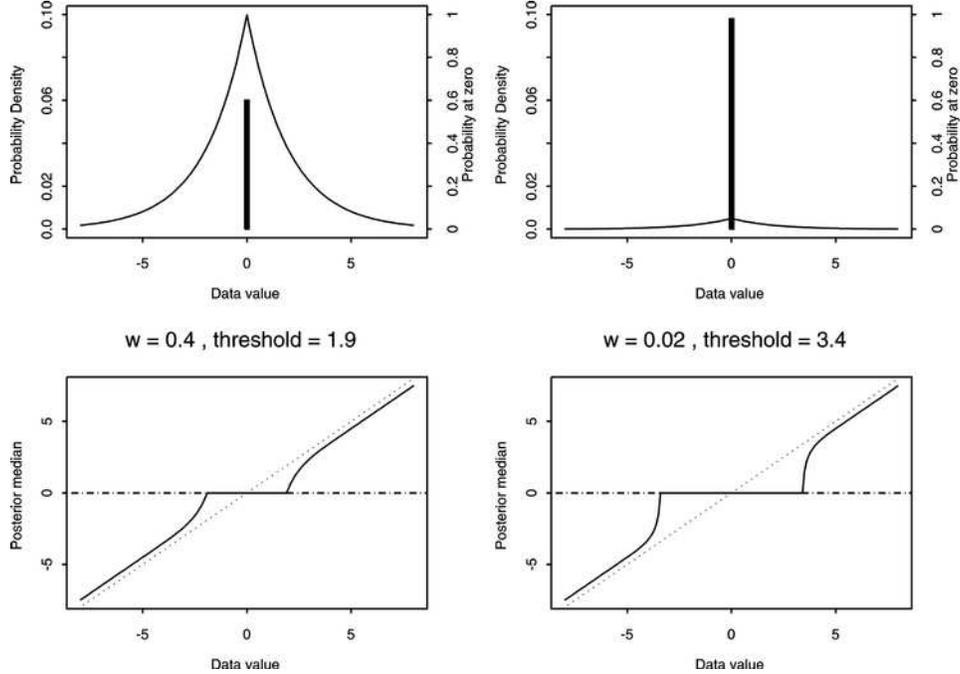

Fig. 5. *First line: Prior distribution for $w = 0.4$ and $w = 0.02$, for the mixed Laplace prior with $a = 0.5$. The atom of probability at zero is represented by the solid vertical bar, plotted to the scale indicated on the right of the plot; the probability density of the nonzero part of the prior is plotted to the scale at the left. Second line: Posterior median functions for the same priors. The dotted line is the diagonal $y = x$. It can be seen that the posterior median is a monotonic function of the data value and is zero whenever the absolute value of the datum is below the threshold.*

subject to the constraint on $w$ that the threshold satisfies $t(w) \leq \sqrt{2 \log n}$. The threshold chosen by the method will then be the value $t(w)$.

The function $\ell'(w)$ is a monotonic function of $w$, so its root is very easily found numerically, provided the function $g$ is tractable; see Section 2.2. Our basic approach will then be to plug the value $\hat{w}$ back into the prior and then estimate the parameters $\mu_i$ using this value of $w$, either using the posterior median itself, or by using some other thresholding rule with the same threshold $t(w)$. In the example above simple hard thresholding was used.

Another possibility is to use the posterior mean, which we denote $\tilde{\mu}(x; w)$, so that the corresponding estimate is $\tilde{\mu}_i = \tilde{\mu}(X_i; \hat{w})$. The posterior mean rule fails to have the thresholding property, and, hence, produces estimates in which, essentially, all the coefficients are nonzero. Nevertheless, it has shrinkage properties that allow it to give good results in certain cases. We shall see that both in theory and in simulation studies, the performance of the posterior mean is good, but not quite as good as the posterior median.

The empirical Bayes is a fully automatic practical method; intuitively, the reason it works well is as follows. If the means $\mu_i$ are all near zero, then $\hat{w}$ will be small, corresponding to a large threshold $t(\hat{w})$, so that most of the means will be estimated to be zero. On the other hand, if the $\mu_i$ are larger, then a small threshold will be chosen, and the data will not be shrunk so severely in the estimation of the vector of means.

1.3. *Measures of sparsity and minimax rates.* The sparsity of a signal is not just a matter of the proportion of $\mu_i$ that are zero or very near zero, but also of more subtle ways in which the energy of the signal $\mu$ is distributed among the various components. Our theory will demonstrate that the empirical Bayes choice of estimated threshold yields a highly adaptive procedure, with excellent properties for a wide range of conditions on the underlying signal.

A natural notion of sparsity is the possibility that $\mu$ is a *nearly black* signal, in the sense that the number of indices $i$ for which $\mu_i$ is nonzero is bounded. We define

$$(3) \qquad \ell_0[\eta] = \left\{ \mu : n^{-1} \sum_{i=1}^n I[\mu_i \neq 0] \leq \eta \right\}.$$

With just the knowledge that $\mu$ falls in $\ell_0[\eta]$, how well can $\mu$ be estimated? Define the minimax average square error by

$$R_{n,2}(\ell_0[\eta]) = \inf_{\hat{\mu}} \sup_{\mu \in \ell_0[\eta]} n^{-1} \sum_{i=1}^n E(\hat{\mu}_i - \mu_i)^2.$$

Donoho, Johnstone, Hoch and Stern (1992) show that, considering $\eta = \eta_n \to 0$ as $n \to \infty$, $R_{n,0}$ is $2\eta(\log \eta^{-1})(1 + o(1))$.

A more subtle characterization of sparsity will not require any $\mu_i$ to be exactly zero, but still constrain most of the energy to be concentrated on a few of the $\mu_i$, by placing bounds on the $p$-norm of $\mu$ for $p > 0$. There are various intuitive ways of understanding why $\|\mu\|_p = (\sum |\mu_i|^p)^{1/p}$ for small $p$ is related to the sparsity of $\mu$. Perhaps the simplest is to consider the energy (the sum of squares) of a vector with $\|\mu\|_p = 1$ for some small $p$. If only one component of $\mu$ is nonzero, then the energy will be 1. If, on the other hand, all the components are equal, then the energy is $n^{1-2/p}$ which tends to zero as $n \to \infty$ if $p < 2$, rapidly if $p$ is near zero. By extension of these examples, if $p$ is small, the only way for a signal in an $\ell_p$ ball with small $p$ to have large energy is for it to consist of a few large components, as opposed to many small components of roughly equal magnitude. Put another way, among all signals with a given energy, the sparse ones are those with small $\ell_p$ norm.

In this case we suppose the signal belongs to an $\ell_p$ norm ball of small radius $\eta$,

$$(4) \qquad \ell_p[\eta] = \left\{ \mu : n^{-1} \sum |\mu_i|^p \leq \eta^p \right\},$$

and define the minimax square error

$$R_{n,2}(\ell_2[\eta]) = \inf_{\hat{\mu}} \sup_{\mu \in \ell_0[\eta]} n^{-1} \sum_{i=1}^{n} E(\hat{\mu}_i - \mu_i)^2.$$

Again, considering $\eta \to 0$ as $n \to \infty$, Donoho and Johnstone (1994) show that, for $p \leq 2$, $R_{n,2}(\ell_p[\eta])$ is $\eta^p (2 \log \eta^{-p})^{1-p/2}(1 + o(1))$.

The estimator that attains the ideal performance over a nearly black class, or over an $\ell_p$ ball for some $p > 0$, will in general depend on $p$ and on $\eta$. The minimax rate is a benchmark for the estimation of signals that display the sparseness characteristic of membership of an $\ell_p$ class. Our main theorem will show that, under mild conditions, an empirical Bayes thresholding estimate will essentially achieve the minimax rate over $\eta$ simultaneously for all $p$ in $[0, 2]$, including the nearly black class as the case $p = 0$. In this sense it adapts automatically to the degree and character of sparsity of the signal in the optimum possible way.

A particular minimax risk is the risk when there is no constraint at all on the underlying signal. In this case the minimax asymptotic risk is a constant 1, for example, achieved by the estimator that simply estimates $\mu_i$ by $X_i$. We show that the maximum possible risk of the empirical Bayes thresholding method, under appropriate conditions, is also uniformly bounded, so the adaptivity is not bought at the price of asymptotically unbounded risk for signals of certain kinds.

### 1.4. *Robustness.*

While adaptivity of an estimator is obviously desirable, it is also important that the estimator should be robust to assumptions made. There are several aspects of such robustness that we demonstrate for the empirical Bayes threshold estimator.

*Assumptions on the signal*: Although our procedure is derived from the sparse prior model (2), we derive results under the much weaker assumption that the underlying signal belongs to an appropriate $\ell_p$ ball.

*Assumptions on the noise*: For example, in Section 5 we relax the assumption of Gaussian errors in order to investigate the relation between tails of the prior and tails of the noise density. While, in their present form, some other aspects of our subsequent discussion make use of Gaussian assumptions, the key properties of the posterior median thresholding rule hold under considerably weaker assumptions.

*Assumptions on the error measure*: Rate-optimal risk bounds are established for mean $q$th power error measures for all $q \in (0, 2]$, not just the standard mean square error. Excessive reliance on mean square error ($q = 2$) is often criticized, for example, as not corresponding to visual assessments of error. Choices of $q < 2$ will give greater (relative) weight to small errors, and in some sense, the $q \to 0$ limit corresponds to counting the number of errors $I\{\hat{\mu}_i \neq \mu_i\}$.

*Assumptions on the estimator itself*: While the posterior median is the motivating estimator for our work, the exact form of the thresholding rule is not specified in our theoretical discussion. The key point is that the data dependent threshold is chosen according to the sparse empirical Bayes prescription. Indeed, the processing rule does not even have to be a strict thresholding rule. We obtain good results for the posterior mean, which is not a thresholding rule but still possesses an appropriate bounded shrinkage property; however, for full robustness to the choice of error measure, strict thresholding rules have to be used.

1.5. *Related work.* Abramovich, Benjamini, Donoho and Johnstone (2000) show that the false discovery rate approach provides adaptive asymptotic minimaxity at the level of exact constants, as well as the rates of convergence that we demonstrate for the empirical Bayes method. However, their results do not guarantee robustness for denser signals, and there is some evidence of this nonrobustness in the simulations we report in Section 3.

In a more restrictive scenario than ours, and mainly concentrating on the application to wavelet smoothing, Zhang (2004) provides an asymptotically more sharply adaptive empirical Bayes analysis. This analysis uses much more general families of priors than our simple mixtures, and employs nonparametric infinite-order kernel methods to estimate the corresponding marginal densities. Such methods are complex to implement in software, and their sharp asymptotic properties might not be apparent in moderate samples.

Mixture priors built from models such as (2) are quite common in Bayesian variable selection problems: our interest was stimulated in part by analysis of a proposal due to George and Foster (1998, 2000) which takes $\gamma$ to be Gaussian. For further references specifically in the wavelet setting, see the companion paper Johnstone and Silverman (2004).

1.6. *Outline of the paper.* The paper now proceeds as follows. In Section 2 we set out some key definitions and state the main theorem of the paper. To show that the advantages of the estimate are not just theoretical, in Section 3 a simulation study is presented, comparing the empirical Bayes method with a range of other estimators, on cases covering both sparse and dense signals. In this study the theoretical adaptivity and robustness properties of the empirical Bayes method are clearly borne out. In very sparse cases the theory suggests that some asymptotic improvement may be possible for very sparse signals, and in Section 4, we set out a modification of our standard procedure, whereby the threshold is increased by a suitable factor when the signal is estimated to be very sparse. We state a result giving key properties of this procedure, and also present some discussion and numerical results that suggest that, except when the sample size is very large indeed, the modification may be of theoretical interest only.

We then move to the proofs of the main results. In Section 5 various detailed preliminaries are considered, including the properties of the posterior rules under more general noise distributions than the Gaussian. We then go on, in Section 6, to consider risk bounds first for fixed thresholds, and then for data-dependent thresholds. These bounds depend on tail probabilities for the random thresholds. As a prerequisite to the control of these probabilities, Section 7 investigates properties of the moment behavior of the marginal likelihood score function. In Section 8 the proof of the main theorem is completed: the results of Section 7 yield tail probabilities of the prior parameters chosen by the empirical Bayes method, and, hence, of the corresponding random thresholds. These are fed into the bounds of Section 6 to complete the proof. Section 9 then contains the modifications to the previous arguments needed to prove Theorem 2.

The conditions in the main theorem for the posterior mean do not cover as wide a range of loss functions as for strict thresholding rules. In Section 10 it is shown that this is an essential feature of the use of such a rule; for values of $q \leq 1$ the posterior mean cannot yield an optimal estimate relative to $q$th power loss under the same broad conditions.

## 2. Aspects of the sequence estimation problem.

It is convenient to set up some notational conventions. Where $A_r$ and $B_r$ are numerical quantities depending on a discrete or continuous index $r$, we write $A_r \asymp B_r$ to denote $0 < \liminf_r A_r/B_r \leq \limsup_r A_r/B_r < \infty$, and $A_r \sim B_r$ to denote $A_r/B_r \to 1$. We use $\phi$ and $\Phi$ for the standard normal density and cumulative, respectively, and set $\bar{\Phi} = 1 - \Phi$. When there is no confusion about the value of the prior weight $w$, it may be suppressed in our notation. Use $c$ and $C$ to denote generic strictly positive constants, not necessarily the same at each use, even within a single equation. We adopt the convention that $c$ is an absolute constant, while the use of $C$ will indicate a possible dependence on the prior density component $\gamma$.

### 2.1. Assumptions on the prior.
When using the mixture prior (2), we shall see that there are considerable advantages in using a heavy-tailed density for $\gamma$, for example, the Laplace density

$$(5) \qquad \gamma(u) = \tfrac{1}{2}\exp(-|u|)$$

or the mixture density given by

$$(6) \qquad (\mu|\Theta = \theta) \sim N(0, \theta^{-1} - 1) \qquad \text{with } \Theta \sim \text{Beta}(\alpha, 1).$$

The latter density for $\mu$ has tails that decay as $\mu^{-2\alpha - 1}$, so that, in particular, if $\alpha = \frac{1}{2}$, then the tails will have the same weight as those of the Cauchy distribution. To be explicit, this has

$$\gamma(u) = \int_0^1 \frac{1}{\sqrt{8\pi(1-\theta)}} \exp\left\{-\frac{1}{2} u^2 \theta (1-\theta)^{-1}\right\} d\theta.$$

In both cases (5) and (6) the posterior distribution of $\mu$ given an observed $X$, and the marginal distribution of $X$, are tractable, so that the choice of $w$ by marginal maximum likelihood, and the estimation of $\mu$ by posterior mean or median, can be performed in practice, using the approach outlined in Section 2.2. Details of the relevant calculations for particular priors are given by Johnstone and Silverman (2004).

Throughout the paper we will assume that the nonzero part of the prior, $\gamma$, has a fixed unimodal symmetric density. In addition, we will assume that

$$(7) \qquad \sup_{u > 0} \left| \frac{d}{du} \log \gamma(u) \right| = \Lambda < \infty.$$

It follows from this assumption that, for $u > 0$, $\log \gamma(u) \geq \log \gamma(0) - \Lambda u$, so that, for all $u$,

$$(8) \qquad \gamma(u) \geq \gamma(0) e^{-\Lambda |u|}.$$

Thus, the tails of $\gamma$ have to be exponential or heavier, and the Gaussian model for $\gamma$ is ruled out. We will also assume that the tails of $\gamma$ are no heavier than Cauchy, in the sense that $u^2 \gamma(u)$ is bounded over all $u$. Finally, we make the mild regularity assumption that, for some $\kappa \in [1, 2]$,

$$(9) \qquad \gamma(y)^{-1} \int_y^\infty \gamma(u)\, du \asymp y^{\kappa - 1} \qquad \text{as } y \to \infty.$$

If $\gamma$ has asymptotically exponential tails, then $\kappa = 1$. If $\gamma(y) \asymp y^{-2}$ for large $y$, then the tail probability is asymptotic to $y^{-1}$ and $\kappa = 2$. Any Pareto tail behavior gives the value $\kappa = 2$.

2.2. *Finding the estimate.* Define the score function $S(w) = \ell'(w)$, and define

$$(10) \qquad \beta(x) = \frac{g(x)}{\phi(x)} - 1 \quad \text{and} \quad \beta(x, w) = \frac{\beta(x)}{1 + w\beta(x)},$$

so that

$$(11) \qquad S(w) = \sum_{i=1}^n \frac{g(X_i) - \phi(X_i)}{(1-w)\phi(X_i) + wg(X_i)} = \sum_{i=1}^n \beta(X_i, w).$$

Since by elementary calculus $\beta(x, w)$ is a decreasing function of $w$ for each $x$, the function $S(w)$ is also decreasing. Let $w_n$ be the weight that satisfies $t(w_n) = \sqrt{2 \log n}$. If $S(w_n) > 0$ and $S(1) < 0$, then the zero of $S$ in the range $[w_n, 1]$ is the estimated weight $\hat{w}$. Furthermore, the sign of $S(w)$ for any particular $w$ specifies on which side of $w$ the estimate $\hat{w}$ lies. [Note that $S$ will be strictly decreasing except in the pathological case where $\beta(X_i) = 0$ for all $i$, when $S(w) = 0$ for all $w$ and the likelihood is constant.]

The marginal maximum likelihood approach can be used to estimate other parameters of the prior. In particular, if a scale parameter $a$ is incorporated by considering a prior density $(1 - w)\delta_0(\mu) + wa\gamma(a\mu)$, define $g_a$ to be the convolution of $a\gamma(a \cdot)$ with the normal density. Then both $a$ and $w$ can be estimated by finding the maximum over both parameters of

$$\ell(w, a) = \sum_{i=1}^{n} \log\{(1 - w)\phi(X_i) + w g_a(X_i)\}.$$

If $\gamma$ is the Laplace density, the tractability of the procedure is not affected by the inclusion of a scale parameter into the prior. In this case if one is maximizing over both $w$ and $a$, then a package numerical maximization routine that uses gradients has been found to be an acceptably efficient way of maximizing $\ell(w, a)$.

In the current paper we will not develop theory for the case where additional parameters of $\gamma$ are estimated, but we will include the possibility of estimating a scale parameter in the simulation study reported in Section 3.

The R/S-PLUS software package EbayesThresh [Johnstone and Silverman (2003)] includes a routine that performs empirical Bayes thresholding on a vector of data. It allows the use of either the Laplace or the quasi-Cauchy prior, and in the case of the Laplace prior, the scale parameter can if desired be chosen by marginal maximum likelihood. Estimation may be carried out using the posterior median or posterior mean rule, or by hard or soft thresholding. In addition, there are several routines that will allow users to develop other aspects of the general approach.

2.3. *Shrinkage rules.* We begin with some definitions, leading up to the statement of the main theorem of the paper. A function $\delta(x, t)$ will be called a *shrinkage rule* if and only if $\delta(\cdot, t)$ is antisymmetric and increasing on $(-\infty, \infty)$ for each $t \geq 0$, and

$$(12) \qquad 0 \leq \delta(x, t) \leq x \qquad \text{for all } x \geq 0.$$

The shrinkage rule $\delta(x, t)$ will be a *thresholding rule* with threshold $t$ if and only if

$$(13) \qquad \delta(x, t) = 0 \qquad \text{if and only if } |x| \leq t,$$

and will have the *bounded shrinkage property relative to the threshold* $t$ if, for some constant $b$,

$$(14) \qquad |x - \delta(x, t)| \le t + b \qquad \text{for all } x \text{ and } t.$$

For any given weight $w$, the posterior median will be a thresholding rule and will have the bounded shrinkage property if $|(\log \gamma)'|$ is bounded; see Lemma 2(v). In Section 5.5 it is demonstrated that the posterior mean for the same weight will have the same bounded shrinkage property, but will not be a strict thresholding rule. If the hyperparameter $w$ is chosen by marginal maximum likelihood, both are examples of rules with random threshold $\hat{t} = t(\hat{w})$.

2.4. *Risk measures and the main result.* As already mentioned, we do not restrict attention to losses based on squared errors, but we measure risk by the average expected $q$th power loss

$$(15) \qquad R_q(\hat{\mu}, \mu) = n^{-1} \sum_{i=1}^{n} E|\hat{\mu}_i - \mu_i|^q, \qquad 0 < q \le 2.$$

Note that the posterior median and mean estimators for prior (2) are Bayes rules for the $q = 1$ and $q = 2$ error measures, respectively.

We set two goals for estimation using the empirical Bayes threshold: "uniform boundedness of risk" and "flexible adaptation." To explain what we mean by flexible adaptation, suppose that the signal is sparse in the sense of belonging to an $\ell_p$ norm ball $\ell_p[\eta]$ as defined in (4). As before, we include nearly black classes as the case $p = 0$. If the radius $\eta$ is small, we would hope that the estimation error $R_q(\hat{\mu}, \mu)$ should be appropriately small. How small is benchmarked in terms of the minimax risk

$$R_{n,q}(\ell_p[\eta]) = \inf_{\hat{\mu}} \sup_{\mu \in \ell_p[\eta]} R_q(\hat{\mu}, \mu).$$

Suppose $\eta = \eta_n \to 0$ as $n \to \infty$ but that, in the case $q > p > 0$,

$$(16) \qquad n^{-1/p} \eta^{-1} (\log \eta^{-p})^{1/2} \to 0,$$

which prevents $\eta$ from becoming very small too quickly. (For $p = 0$ we require $n\eta \to \infty$.) Then we have the asymptotic relation

$$(17) \qquad R_{n,q}(\ell_p[\eta_n]) \sim r_{p,q}(\eta_n) \qquad \text{as } n \to \infty,$$

where

$$(18) \qquad r_{p,q}(\eta) = \begin{cases} \eta^q, & 0 < q \le p, \\ \eta^p (2 \log \eta^{-p})^{(q-p)/2}, & 0 < p < q, \\ \eta (2 \log \eta^{-1})^{q/2}, & p = 0, \ q > 0. \end{cases}$$

The relation (17) is proved by Donoho and Johnstone (1994) for the case $p > 0$ and $q \geq 1$, but only minor modifications are needed to extend the result to all the cases we consider.

We can now state our main result, which gives comparable bounds on the risk function of the empirical Bayes thresholding procedure. Apart from an error of order $n^{-1}(\log n)^{2+(q-p)/2}$, the procedure uniformly attains the same error rate as the minimax estimator for all $p$ in $[0, 2]$ and $q$ in $(0, 2]$.

THEOREM 1. *Suppose that $X \sim N_n(\mu, I)$, that $\delta(x, t)$ is a thresholding rule with the bounded shrinkage property and that $0 \leq p \leq 2$ and $0 < q \leq 2$. Let $\hat{w}$ be the weight chosen by marginal maximum likelihood for a mixture prior* (2) *with $\gamma$ satisfying the assumptions set out in Section* 2.1. *Let $\hat{t} = t(\hat{w})$, where $t(w)$ denotes the threshold of the posterior median rule corresponding to the prior weight $w$. Then the estimator $\hat{\mu}_i(x) = \delta(x_i, t(\hat{w}))$ satisfies:*

(a) *(Uniformly bounded risk) There exists a constant $C_0(q, \gamma)$ such that*

$$\sup_{\mu} R_q(\hat{\mu}, \mu) \leq C_0.$$

(b) *(Adaptivity) There exist constants $C_i(p, q, \gamma)$ such that for $\eta \leq \eta_0(p, q, \gamma)$ and $n \geq n_0(p, q, \gamma)$,*

$$\sup_{\mu \in \ell_p[\eta]} R_q(\hat{\mu}, \mu) \leq C_1 r_{p,q}(\eta) + C_2 n^{-1}(\log n)^{2+(q-p)/2}. \tag{19}$$

*When $q \in (1, 2]$, these results also hold for the posterior mean estimate $\tilde{\mu}$.*

We emphasize that it is not necessary that $\delta(x, t)$ be derived from the posterior median or mean rule. It might be hard or soft thresholding or some other nonlinearity with the stated properties. The point of the theorem is that empirical Bayes estimation of the threshold parameter suffices with all such methods to achieve both adaptivity and uniformly bounded risk.

If $q > p > 0$, then we necessarily have $p < 2$, and the first term of (19) dominates if $\eta^p > n^{-1} \log^2 n$ and the second if $\eta^p < n^{-1} \log^2 n$. It follows that the result is equivalent to

$$\sup_{\mu \in \ell_p[\eta]} R_q(\hat{\mu}, \mu) \leq \begin{cases} C r_{p,q}(\eta), & \text{if } \eta^p \geq n^{-1} \log^2 n, \\ C n^{-1}(\log n)^{2+(q-p)/2}, & \text{if } \eta^p < n^{-1} \log^2 n. \end{cases} \tag{20}$$

Note that $\eta^p \geq n^{-1} \log^2 n$ is a sufficient condition for (16). For the nearly black case $p = 0$, a similar argument leads to (20) with $\eta^p$ replaced by $\eta$.

If $p \geq q$, the bound can be written as

$$\sup_{\mu \in \ell_p[\eta]} R_q(\hat{\mu}, \mu) \leq C \max\{\eta^q, n^{-1}(\log n)^{2+(q-p)/2}\} \tag{21}$$

and the "break-even" point between the two bounds occurs at a value of $\eta$ bounded above by $\eta^p = n^{-1} \log^2 n$. It remains the case that for $\eta^p \geq n^{-1} \log^2 n$ the supremum of the risk is bounded by a multiple of $r_{p,q}(\eta)$. Therefore, for every $p$ and $q$ in $(0, 2]$, and for the nearly black case $p = 0$, our estimator attains the optimal $q$-norm risk (18), up to a constant multiplier, for all sufficiently large $n$ and for $\eta$ satisfying $n^{-1} \log^2 n \leq \eta^p \leq \eta_0^p$ if $p > 0$ and $n^{-1} \log^2 n \leq \eta \leq \eta_0$ if $p = 0$.

**3. Some simulation results.** In order to investigate the capability of the empirical Bayes method to adapt to the degree of sparsity in the true signal, a simulation study was carried out. We approach the issue of sparsity directly, by explicitly constructing sequences with a wide range of sparse behavior. The S-PLUS code used to carry out the simulations is available from the authors' web sites, enabling the reader both to verify the results and to conduct further experiments if desired.

As an initial range of models for sparse behavior, we fixed the sample size $n$ to 1000. We considered the estimation of a sequence $\mu$ which has $\mu_i = 0$ except in $K$ randomly chosen positions, where it takes a specified value $\mu_0$. For each $i$, a data value $X_i \sim N(\mu_i, 1)$ is generated, and various methods are used to estimate the sequence $\mu$ from the sequence of $X_i$.

The parameter $K$ controls the sparsity of the signal, and the values for which results are reported are 5, 50 and 500—ranging from a very sparse signal, indeed, to one in which half the data contain nonzero signal. The other parameter $\mu_0$ gives the strength of the signal if it is nonzero. The values reported were 3, 4, 5 and 7, bearing in mind that the noise is $N(0,1)$. One hundred replications were carried out for each of the values of $K$ and $\mu_0$, with the same 100,000 noise variables used for each set of replications.

The posterior median estimator was used, with the prior parameters chosen by marginal maximum likelihood for two different functions $\gamma$ for the nonzero part of the prior. The double exponential $\gamma(u) = \frac{1}{2} a \exp(-a|u|)$ was used with both the parameter $a$ and the prior weight $w$ chosen by marginal maximum likelihood. For comparison, the heavy-tailed mixture density with Cauchy tails, as defined in (6) with $\alpha = \frac{1}{2}$, was also considered. For both choices of the function $\gamma$, the performance of the posterior median as a point estimator was studied. For double exponential $\gamma$ with both parameters estimated, two other estimators were also considered, the posterior mean, and hard thresholding with threshold equal to that of the posterior median function. In addition, the effect of fixing the scale parameter in the double exponential was investigated by considering four different values of $a$; in each case $w$ was chosen by marginal maximum likelihood and the posterior median estimator used.

These methods were compared with classical soft and hard universal thresholding (using the threshold $\sqrt{2 \log n} \approx 3.716$) and with three other methods intended to be adaptive to different levels of sparsity.

The SURE method [Donoho and Johnstone (1995)] aims to minimize the mean squared error of reconstruction, by minimizing Stein's unbiased risk estimate for the mean squared error of soft thresholding. Thus, we choose $\hat{t}_{\text{SURE}}$ as the minimizer (within the range $[0, \sqrt{2 \log n}]$) of

$$\hat{U}(t) = n + \sum_1^n x_k^2 \wedge t^2 - 2 \sum_1^n I\{x_k^2 \leq t^2\}.$$

This is based on the unbiased risk estimator of Stein (1981) in the estimation of a multivariate normal mean. In addition, a modification proposed by Donoho and Johnstone (1995) aimed at gaining greater adaptivity is considered; this chooses between the SURE and universal thresholds according to the result of a test for sparsity; see also Section 6.4.2 of Bruce and Gao (1996) for details.

The false discovery rate (FDR) approach is derived from the principle of controlling the false discovery rate in simultaneous hypothesis testing [Benjamini and Hochberg (1995)] and has been studied in detail in the estimation setting, for example, by Abramovich, Benjamini, Donoho and Johnstone (2000). Order the data by decreasing magnitudes:

$$|x|_{(1)} \geq |x|_{(2)} \geq \cdots \geq |x|_{(n)}$$

and compare to a *quantile boundary*:

$$t_k = \sigma z(q/2 \cdot k/n),$$

where the false discovery rate parameter $q \in (0, 1/2)$. Define a crossing index

$$\hat{k}_F = \max \{k : |x|_{(k)} \geq t_k\}$$

and use this to set the threshold $\hat{t}_F = t_{\hat{k}_F}$. Various values for the rate parameter $q$ were used.

Block thresholding methods are designed to make use of neighboring information in setting the threshold applied to each individual data point. We considered the BlockThresh method of Cai (2002) and the hard thresholding versions of the NeighBlock and NeighCoeff methods of Cai and Silverman (2001). The principle of all these methods is to consider the data in blocks. BlockThresh thresholds all the data in each block by reference to the sum of squares of the data in the block. The other two methods use overlapping blocks and keep or zero the data in the middle of each block according to the sum of squares over the whole block. See the original papers for more details.

For each method considered, for each replication the total squared error of the estimation $\sum (\hat{\mu}_i - \mu_i)^2$ was recorded, and the average over 100 replications is reported. The square error of every replication is available from

TABLE 1

*Average of total squared error of estimation of various methods on a mixed signal of length 1000. A given number of the original signal values is set equal to a nonzero value, and the remainder are zero. In each column those entries that outperform the MML/exponential/posterior median method are underlined. Those that outperform by more than about 10% are set in bold type. The row marked "postmean" refers to the posterior mean using the double exponential model. The row "exphard" refers to hard thresholding using the threshold given by the posterior median of the marginal maximum likelihood choice within the double exponential model. The rows for fixed values of a correspond to the posterior median where only the weight w is chosen by MML and the scale parameter a is fixed at the given value*

| Number nonzero | 5 | | | | 50 | | | | 500 | | | |
|---|---|---|---|---|---|---|---|---|---|---|---|---|
| Value nonzero | 3 | 4 | 5 | 7 | 3 | 4 | 5 | 7 | 3 | 4 | 5 | 7 |
| Exponential | 36 | 32 | 17 | 8 | 214 | 156 | 101 | 73 | 857 | 873 | 783 | 658 |
| Cauchy | 37 | 36 | 18 | 8 | 271 | 176 | 103 | 77 | 922 | 898 | 829 | 743 |
| Postmean | 34 | 32 | 21 | 11 | 201 | 169 | 122 | 85 | 860 | 888 | 826 | 708 |
| Exphard | 51 | 43 | 22 | 11 | 273 | 189 | 130 | 91 | 998 | 998 | 983 | 817 |
| $a = 1$ | 36 | 32 | 19 | 15 | 213 | 166 | 142 | 135 | 994 | 1099 | 1126 | 1130 |
| $a = 0.5$ | 37 | 34 | 17 | 10 | 244 | 158 | 105 | 92 | 845 | 878 | 884 | 884 |
| $a = 0.2$ | 38 | 37 | 18 | 7 | 299 | 188 | 95 | 69 | 1061 | **730** | **665** | 656 |
| $a = 0.1$ | 38 | 37 | 18 | **6** | 339 | 227 | 102 | **60** | 1496 | 798 | **609** | **579** |
| SURE | 38 | 42 | 42 | 43 | 202 | 209 | 210 | 210 | 829 | 835 | 835 | 835 |
| Adapt | 42 | 63 | 73 | 76 | 417 | 620 | 210 | 210 | 829 | 835 | 835 | 835 |
| FDR $q = 0.01$ | 43 | 51 | 26 | **5** | 392 | 299 | 125 | **55** | 2568 | 1332 | **656** | **524** |
| FDR $q = 0.1$ | 40 | 35 | 19 | 13 | 280 | 175 | 113 | 102 | 1149 | **744** | **651** | 644 |
| FDR $q = 0.4$ | 58 | 58 | 53 | 52 | 298 | 265 | 256 | 254 | 919 | 866 | 860 | 860 |
| BlockThresh | 46 | 72 | 72 | 31 | 444 | 635 | 600 | 293 | 1918 | 1276 | 1065 | 983 |
| NeighBlock | 47 | 64 | 51 | 26 | 427 | 543 | 439 | 227 | 1870 | 1384 | 1148 | 972 |
| NeighCoeff | 55 | 51 | 38 | 32 | 375 | 343 | 219 | 156 | 1890 | 1410 | 1032 | 870 |
| Universal soft | 42 | 63 | 73 | 76 | 417 | 620 | 720 | 746 | 4156 | 6168 | 7157 | 7413 |
| Universal hard | 39 | 37 | 18 | **7** | 370 | 340 | 163 | **52** | 3672 | 3355 | 1578 | **505** |

the authors' web sites for any reader who wishes to examine the results in more detail.

Some results are given in Table 1 and the following conclusions can be drawn:

- The Cauchy method is always nearly, but not quite, as good as the exponential method. Our theory is not sensitive enough to discriminate between the two methods.

- In general, the posterior mean does not perform quite as well as the posterior median.
- It is better to use the posterior median function itself rather than hard thresholding with the resulting threshold.
- In the case $\mu_0 = 7$ where the nonzero signal is very clearly different from zero, hard thresholding with the universal threshold performs somewhat better than the exponential method, but in other cases, particularly with moderate or large amounts of moderate sized signal, it can give disastrous results.
- Estimating the scale parameter $a$ is probably preferable to using a fixed value, though it does lead to slower computations. In general, the automatic choice is quite good at tracking the best fixed choice, especially for a sparse and weak signal.
- SURE is a competitor when the signal size is small ($\mu_0 = 3$), but performs poorly when $\mu_0$ is larger, particularly in the sparser cases. The attempt to make SURE more adaptive is counterproductive.
- If $q$ is chosen appropriately, FDR can outperform exponential in some cases, but the choice of $q$ is crucial and varies from case to case. With the wrong choice of $q$, the performance of FDR can be poor.
- The block thresholding methods do not perform very well. In the companion paper [Johnstone and Silverman (2004)] block thresholding methods are also compared with empirical Bayes methods for the thresholding of wavelet coefficients, and the difference in performance is not so great. This is presumably because there is some correlation among the positions in which the wavelet coefficients are effectively nonzero. By contrast, in the test signals under current consideration, the nonzero positions are chosen by uniform random sampling without replacement.
- The median standard error of the entries of the table with 5 nonzero coefficients is around 1, with corresponding figures of about 3 for those with 50 nonzero coefficients, and 5 for the entries with 500 nonzero coefficients. Generally speaking, the standard errors tend to be smaller for the empirical Bayes methods than for the other methods considered; the false discovery rate and block thresholding methods have errors that have variance two to three times as large as the double exponential MML posterior median method, and for the universal thresholding methods the variance is higher by a factor of about 5. This is an indication of the stability of the empirical Bayes methods.
- Not surprisingly, given that the same data are used for all cases, the standard error of the comparison between the first method and the other methods in the table is typically smaller than that for individual entries taken alone. The comparison standard error has a median value of 0.8 for the sparsest signals and about 2 for the signals with 50 and 500 nonzero elements. In general, comparisons between empirical Bayes methods have somewhat smaller standard errors than those involving other approaches. Only about 10% of the comparisons between the top line and other entries in the table are within 3 standard errors of zero, and all the comparisons that are numerically more than trivial are clearly statistically significant on the basis of the study we have carried out.

The two SURE methods, the FDR method with $q = 0.01$ or $q = 0.4$, and the two universal thresholding methods all have the property that there is a case in which their measured error is around three or more times that of the exponential method, while never, or hardly ever offering any substantial improvement. Hence, all are much worse at adapting to different patterns of sparsity. The FDR method with $q = 0.1$ is a better competitor, but only wins in four of the twelve cases. The best improvement over exponential is 651/783, a 17% improvement, while the best improvement of the exponential over the FDR method is 73/102, nearly 30%. Taking both adaptivity and overall performance into account, the exponential is clearly the estimator of choice.

In order to quantify the comparison between the various methods, for each of the models considered define the *inefficiency* of a method A for a particular model B to be

$$100 \times \left[ \frac{\text{average error for method A applied to data from model B}}{\text{minimum error for any method for model B}} - 1 \right].$$

Twelve different models are considered in Table 1, and summary statistics for the twelve inefficiency values for the various methods are given in Table 2. The posterior median of the exponential model with estimated scale parameter is the best on nearly every measure: the maximum inefficiency of the Cauchy and exponential ($a = 0.2$) methods is slightly smaller, but both of these methods are decisively beaten on the median inefficiency and are also equaled or beaten on the other two measures.

**4. Modifying the threshold for very sparse signal.** In this section we discuss a possible modification of the estimator, which allows a reduction in error in very sparse cases, when the overwhelming majority of components have essentially zero signal. Our original motivation for this arises from the use of wavelet methods to estimate derivatives, where it was shown by Abramovich and Silverman (1998) that the appropriate universal threshold is not $\sqrt{2 \log n}$, but is a multiple of this quantity. The basic notion of the modified estimator is this: if the threshold $\hat{t} = t(\hat{w})$ estimated by the marginal maximum likelihood method is at or near the universal threshold, we replace it by a higher threshold.

TABLE 2

*Comparison of methods: for each method the stated median, mean, maximum and tenth largest inefficiency is over the twelve cases considered in Table 1*

|  | median | mean | 10th | max |
|---|---|---|---|---|
| Exponential | 7 | 17 | 30 | 52 |
| Cauchy | 19 | 25 | 42 | 47 |
| Postmean | 22 | 27 | 40 | 95 |
| Exphard | 37 | 46 | 62 | 93 |
| $a = 1$ | 35 | 57 | 124 | 165 |
| $a = 0.5$ | 15 | 29 | 75 | 84 |
| $a = 0.2$ | 18 | 19 | 30 | 48 |
| $a = 0.1$ | 14 | 24 | 45 | 80 |
| SURE | 35 | 121 | 151 | 676 |
| Adapt | 103 | 223 | 303 | 1282 |
| FDR $q = 0.01$ | 44 | 56 | 91 | 210 |
| FDR $q = 0.1$ | 18 | 35 | 39 | 139 |
| FDR $q = 0.4$ | 71 | 169 | 214 | 847 |
| BlockThresh | 129 | 228 | 456 | 531 |
| NeighBlock | 119 | 181 | 335 | 376 |
| NeighCoeff | 106 | 136 | 131 | 486 |
| Universal soft | 529 | 643 | 1282 | 1367 |
| Universal hard | 50 | 100 | 159 | 359 |

4.1. *Definition of the modified estimator and theoretical discussion.* To be precise, set $t_n^2 = 2 \log n - 5 \log \log n$. Let $A \geq 0$ be fixed and put $t_A = \sqrt{2(1 + A) \log n}$. Then define

$$\hat{t}_A = \begin{cases} \hat{t}, & \text{if } \hat{t} \leq t_n, \\ t_A, & \text{if } \hat{t} > t_n. \end{cases}$$

With this modification of the threshold, we can reduce the order of magnitude of the part of the error in Theorem 1, as follows.

THEOREM 2. *Under the assumptions of Theorem 1, let $\hat{\mu}_A$ be defined by $\hat{\mu}_{A,i}(x) = \delta(x_i, \hat{t}_A)$. Define $\mathcal{R}_{p,q}^A(\eta) = \sup_{\mu \in \ell_p[\eta]} R_q(\hat{\mu}_{A,i}, \mu)$. Then, for suitable constants $C$,*

$$(22) \qquad \mathcal{R}_{p,q}^A(\eta) \leq C \qquad \text{for all } \eta$$

*and, for all sufficiently large n and for suitable* $\eta_0$,

$$(23) \quad \mathcal{R}_{p,q}^A(\eta) \leq C \max \{r_{p,q}(\eta), n^{-1-A}(\log n)^{(q-1)/2}\} \qquad \text{for } \eta \leq \eta_0.$$

*For* $q > p > 0$ *we also have, for sufficiently large n,*

$$\mathcal{R}_{p,q}^A(\eta) \leq C \max \{n^{(q-p)/p}\eta^q, n^{-1-A}(\log n)^{(q-1)/2}\}$$

$$(24)$$

$$\text{for } \eta^p < n^{-1}(\log n)^{p/2}.$$

The ramifications of this theorem in the wavelet context are explored by Johnstone and Silverman (2004), but it has independent interest in exposing the different regimes for adaptive estimation, especially in the case $q > p$. Note first that conclusion (22) is the same as for the unmodified estimator, and in the range $\eta^p \geq n^{-1}\log^2 n$ for $p > 0$ ($\eta \geq n^{-1}\log^2 n$ for $p = 0$) so is (23), because in that range the dominating term in the error is $r_{p,q}(\eta)$ for both estimators.

For $q > p > 0$, define $\eta_1^p = n^{-1}(\log n)^{p/2}$. For $\eta > \eta_1$, $r_{p,q}(\eta)$ is bounded by a multiple of $n^{(q-p)/p}\eta^q$ and so (24), in fact, holds for all $\eta < \eta_0$, but only gives a stronger result than (23) if $\eta < \eta_1$. This is not in contradiction with the result (17) of Donoho and Johnstone (1994) because the condition (16) can be rewritten precisely as

$$r_{p,q}(\eta) = o(n^{(q-p)/p}\eta^q),$$

or equivalently, $\eta/\eta_1 \to \infty$.

The three bounds in Theorem 2 may be considered as corresponding to three different zones of estimation. If $\eta > \eta_0$, then the signal is insufficiently sparse for any order of magnitude advantage to be gained by the use of our thresholding method. In the zone $\eta_1 \leq \eta \leq \eta_0$, a suitable thresholding method allows for considerable improvement over the use of a "classical" estimator. Finally, in the extremely sparse zone $\eta < \eta_1$, the $\eta$-dependent part of the error achieved by our estimator compares to that given by the estimator that simply returns the value zero.

Note finally that for the standard estimator all the $\eta$-dependent risks in the zone $\eta^p < n^{-1}\log^2 n$ (for $p > 0$) are dominated by the term $n^{-1}(\log n)^{2+(q-p)/2}$ in the error bound. Since $n^{-1}\log^2 n > \eta_1^p$, the zone where $n^{-1}(\log n)^{2+(q-p)/2}$ dominates includes the very sparse zone, and so there is no point in pursuing the kind of adaptivity within the very sparse zone achieved asymptotically by the modified estimator.

4.2. *Practicalities.* Both theory and intuition suggest that the modified estimator may only be advantageous for very large values of $n$, where $5 \log \log n$ is small relative to $2 \log n$. Otherwise, data that ought to be thresholded with moderate thresholds will essentially be zeroed instead. For example, for $n = 10^6$, we have $5 \log \log n = 13.13$ and $2 \log n = 27.63$. Hence, if the squared estimated threshold in the standard estimator is any more than about half the universal threshold, the modification will use a much larger threshold, thereby causing problems for signals that are nowhere near the very sparse zone.

A version of the modified estimator was investigated by simulation on the same models as considered in Table 1. The Laplace prior with both parameters estimated by marginal maximum likelihood was used. If the estimated threshold was less than 95% of the universal threshold, the posterior median estimate was used. Otherwise, we used hard thresholding with threshold $2\sqrt{\log 1000}$, corresponding to $A = 1$.

The only models for which the estimates were affected were those with only 5 nonzero entries. In each case the average squared error was *increased* by the use of the modified estimator, respectively to 41, 40, 26 and 13, as compared to 36, 32, 17 and 8, for the cases where the nonzero parameter value was 3, 4, 5 and 7. Reducing the number of nonzero parameters to 1 did not change the relative performance of the unmodified and modified methods, unless the nonzero parameter value was also increased. The only case tested where the modified method improved the performance was where there was a single nonzero parameter value with value 10. In this case the unmodified estimator has an average squared error (over 100 simulations) of about 2.4, while the modified estimator has a mean squared error of just over 1. As might be expected, the modified estimator is only clearly advantageous in very sparse cases where nonzero values of the parameters are well above the universal threshold—and in these cases the error of the unmodified method is already very small, so any improvement may be large in relative terms but small in absolute terms.

**5. Proofs of results: some detailed preliminaries.** The remainder of the paper is devoted to the proofs of the theorems stated above. We begin in this section with a detailed discussion of a number of topics that will be useful later in the proof. In some cases these also cast a broader light on the empirical Bayes thresholding procedure. Our proofs cover the cases of nearly black and strong $\ell_p$ constraints on the underlying parameter vector $\mu$. We conjecture that similar arguments can be used for weak $\ell_p$ constraints too, but the full details are left for future investigation.

5.1. *Properties and definitions for the mixture model.* The arguments of this section and the next do not strongly depend on the precise assumption of Gaussian errors in model (1). Indeed, relaxing this assumption sheds some light on the robustness of our results to model formulation; see Remark 1.

For the moment then, we assume that in model (1) the noise coefficients $\epsilon_i$ are i.i.d. from a symmetric Polya frequency $PF_3$ density $\varphi$. Polya frequency functions are discussed in detail by Karlin (1968), and from a statistical perspective by Brown, Johnstone and MacGibbon (1981). The defining property of a $PF_3$ function $\varphi$ is that for $y_1 < y_2 < y_3$ and $z_1 < z_2 < z_3$,

$$(25) \qquad \det_{1 \le i,j \le 3}[\varphi(y_i - z_j)] \ge 0.$$

Examples of such densities include the Gaussian density $\phi$ (observe the distinction of notation), as well as the somewhat heavier tailed Laplace density $\frac{1}{2} e^{-|x|}$ and logistic density $e^{-x}/(1 + e^{-x})^2$. The $PF_3$ assumption implies that $\varphi$ is log-concave, and hence there exists $\rho > 0$ such that

$$(26) \qquad \varphi(y)e^{\rho y} \text{ is decreasing for sufficiently large } y.$$

Thus, the tails of $\varphi$ cannot be heavier than exponential.

For this section only, we also modify assumption (7) on the prior to require only that $\gamma(u) > 0$ for all $u$ and the existence of positive $\Lambda$ and $M$ such that

$$(27) \qquad \sup_{u > M} \left| \frac{d}{du} \log \gamma(u) \right| \le \Lambda < \rho.$$

[In the Gaussian error case, $\varphi = \phi$, this places no essential constraint on $\Lambda$, because we can choose $\rho$ to be arbitrarily large.] Assumptions (26) and (27) taken together imply that the tails of the prior $\gamma$ are heavier than those of the noise density.

The first part of the following lemma shows that the convolution $\gamma \star \varphi$ inherits properties assumed of $\gamma$.

LEMMA 1.  *Assume* (26) *and* (27), *and let* $g = \gamma \star \varphi$. *Then*

$$(28) \qquad g(x) \asymp \gamma(x),$$

$$(29) \qquad (1 + u^2)g(u) \text{ is bounded for all } u,$$

$$(30) \qquad g(y)^{-1} \int_y^\infty g(u)\,du \asymp y^{\kappa - 1},$$

$$(31) \qquad \limsup_{u \to \infty} |(\log g)'(u)| \le \Lambda$$

*and* $g/\varphi$ *is strictly increasing from* $(g/\varphi)(0) < 1$ *to* $+\infty$ *as* $x \to \infty$.

PROOF.  It follows from (27) that $e^{\Lambda y}\gamma(y)$ is an increasing function of $y$ for $y > M$, and since $\gamma$ is unimodal, that for all $x$ and $y$ in $[0, M]$,

$$e^{\Lambda x}\gamma(x) \le C e^{\Lambda y}\gamma(y)$$

for some $C > 1$. Combining these two observations implies that, given any $x > 0$ and $u > 0$,

$$(32) \qquad \gamma(x + u) \ge C^{-1} e^{-\Lambda u} \gamma(x) \quad \text{and} \quad \gamma(x - u) \le C e^{\Lambda u} \gamma(x).$$

It follows that $g(x) \asymp \gamma(x)$, since

$$g(x) > \int_0^\infty \varphi(u)\gamma(x+u)\,du \geq C^{-1}\gamma(x)\int_0^\infty e^{-\Lambda u}\varphi(u)\,du$$

and

$$g(x) = \int_0^\infty \varphi(u)\{\gamma(x+u)+\gamma(x-u)\}\,du \leq C\gamma(x)\int_0^\infty (1+e^{\Lambda u})\varphi(u)\,du,$$

and the right-hand integrals are finite from (26) because $\Lambda < \rho$. Properties (29) and (30) follow immediately from (28) and the assumptions on $\gamma$.

For (31), setting $\Lambda_\infty = \sup |(\log \gamma)'|$, we have $|\gamma'(u)| \leq \Lambda\gamma(u)$ for $u > M$ and $|\gamma'(u)| \leq \Lambda_\infty \gamma(u)$ for $u < M$. Therefore,

$$|(\log g)'(x)| = |g'(x)|/g(x) = \left|\int_{-\infty}^\infty \varphi(x-u)\gamma'(u)\,du\right|/g(x)$$

$$\leq \Lambda \int_M^\infty \varphi(x-u)\gamma(u)\,du/g(x) + \Lambda_\infty \int_{-\infty}^M \varphi(x-u)\gamma(u)\,du/g(x)$$

$$\leq \Lambda + \Lambda_\infty \rho(x),$$

where from (32)

$$\rho(x) = \int_{x-M}^\infty \gamma(x-v)\varphi(v)\,dv/g(x)$$

$$\leq C[\gamma(x)/g(x)]\int_{x-M}^\infty e^{\Lambda v}\varphi(v)\,dv \to 0 \qquad \text{as } x \to \infty.$$

To demonstrate that $g(x)/\varphi(x)$ is increasing on $[0,\infty)$, let $r_u(x) = [\varphi(x+u)+\varphi(x-u)]/\varphi(x)$. Using the symmetry of $\gamma$, we have the representation

$$\frac{g(x)}{\varphi(x)} = \int_0^\infty r_u(x)\gamma(u)\,du, \tag{33}$$

and so it will suffice to show that, for each $u > 0$, $r_u(x)$ is an increasing function of $x$ on $[0,\infty)$. Suppose that $x_2 > x_1 \geq 0$ and consider the defining inequality (25), with $\{y_i\} = \{-x_1, x_1, x_2\}$ and $\{z_i\} = \{-u, 0, u\}$. Subtracting the second row in the determinant from the first and exploiting symmetry of $\phi$ gives

$$0 \leq [\varphi(-x_1+u) - \varphi(x_1+u)]\begin{vmatrix} 1 & 0 & -1 \\ \varphi(x_1+u) & \varphi(x_1) & \varphi(x_1-u) \\ \varphi(x_2+u) & \varphi(x_2) & \varphi(x_2-u) \end{vmatrix}$$

$$= \varphi(x_1)\varphi(x_2)[\varphi(-x_1+u) - \varphi(x_1+u)][r_u(x_2) - r_u(x_1)].$$

Since $\varphi > 0$ and $\varphi(x_1+u) < \varphi(-x_1+u)$, this implies that $r_u(x_2) \geq r_u(x_1)$, as required.

Finally, to show that $g(x)/\phi(x) \to \infty$ as $x \to \infty$, we have, for any $x > 1$, using the result (8),

$$g(x) \geq \int_0^1 \gamma(x-v)\phi(v)\,dv \geq \gamma(x)\int_0^1 \phi(v)\,dv \geq Ce^{-\Lambda|x|}. \qquad \square$$

**Posterior odds.** Write $\text{Odds}(A|B)$ for $P(A|B)/[1-P(A|B)]$. Given $w$, define $\tilde{w}(x,w)$ to be the posterior weight $P(\mu \neq 0|X=x)$, so that the posterior odds $\text{Odds}(\mu \neq 0|x)$ are given by

$$\Omega(x) = \Omega(x,w) = \frac{\tilde{w}(x,w)}{1-\tilde{w}(x,w)} = \frac{w}{1-w}\frac{g}{\varphi}(x).$$

Define $w_0 = (\varphi/g)(0)/[1+(\varphi/g)(0)]$, so that $\Omega(0,w_0) = 1$. For fixed $w$, Lemma 1 shows that $\Omega(x)$ increases from 1 to $\infty$, so that if $w < w_0$, then there exists $\tau(w) > 0$ for which $\Omega(\tau(w),w) = 1$. If we define $\tau(w) = 0$ for $w \geq w_0$, it follows that $w \to \tau(w)$ is a continuous decreasing function of $w \in [0,1]$. We will repeatedly use the function $\tau$ in our subsequent argument.

A simple consequence of these definitions is that for $w_1 < w_0$,

$$(34) \qquad \Omega(\tau(w_1),w) = \frac{w}{1-w}\frac{1-w_1}{w_1} \geq 1, \qquad \text{if } w \geq w_1.$$

Finally, for $x > \tau$, we clearly have

$$(35) \qquad \Omega(x) = \Omega(\tau)\exp\int_\tau^x \{(\log g)' - (\log \varphi)'\}.$$

### 5.2. *Properties of the posterior median.*

LEMMA 2.  *Assume* (26) *and* (27). *The posterior median $\hat{\mu}(x;w)$ is*

  (i) *monotone in $x$: if $x_1 \leq x_2$, then $\hat{\mu}(x_1) \leq \hat{\mu}(x_2)$,*
  (ii) *antisymmetric: $\hat{\mu}(-x) = -\hat{\mu}(x)$,*
  (iii) *a shrinkage rule: $0 \leq \hat{\mu}(x) \leq x$ for $x \geq 0$,*
  (iv) *a threshold rule: there exists $t(w) > 0$ such that $\hat{\mu}(x) = 0$ if and only if $|x| \leq t(w)$,*
  (v) *bounded shrinkage: there exists a constant $b$ such that for all $w,x$,*

$$|\hat{\mu}(x;w) - x| \leq t(w) + b.$$

REMARK 1.  The lemma demonstrates that the posterior median has all the properties needed for the estimation error bounds that will be derived for Gaussian errors in the subsequent sections. The bounded shrinkage property essentially means that rare large observations are more or less reliably assigned to a sparse signal rather than noise in our Bayesian model; conditions

(26) and (27) indicate that a sufficient condition for this is that the tails of the prior be heavier than the tails of the noise distribution. At least in this situation, one may expect the qualitative features of our theory to remain true; it is left to future work to investigate whether there are differences in quantitative thresholds and, perhaps, in rates of convergence.

PROOF OF LEMMA 2. Suppose that $\mu$ has general prior density $f$, with respect to a suitable dominating measure. Then the posterior density

$$f(\mu|x) = C(x)\varphi(x-\mu)f(\mu),$$

so that, for any $u < v$ and $x_2 > x_1$,

$$\frac{f(v|x_2)f(u|x_1)}{f(u|x_2)f(v|x_1)} = \frac{\varphi(x_2-v)}{\varphi(x_2-u)}\frac{\varphi(x_1-u)}{\varphi(x_1-v)} \geq 1,$$

so that

$$f(v|x_1)f(u|x_2) \leq f(u|x_1)f(v|x_2).$$

Now, for any $m$, integrate with respect to the dominating measure over $-\infty < u \leq m$ and $m < v < \infty$ to obtain

$$P(\mu > m|x_1)P(\mu \leq m|x_2) \leq P(\mu > m|x_2)P(\mu \leq m|x_1),$$

so that the odds that $\mu \leq m$ are greater for $X = x_1$ than for $X = x_2$. Letting $\hat{\mu}(x)$ be the posterior median of $\mu$, given $X = x$, it follows that $\hat{\mu}(x_2) \geq \hat{\mu}(x_1)$, so the posterior median is a monotonic function of $x$.

Return now to the mixture prior (2). The antisymmetry of the posterior median is immediate from the symmetry of the prior and the error distribution. If $w > 0$, the probabilities $P(\mu < 0|X = x)$ and $P(\mu = 0|X = x)$ will be nonzero for all $x$ and each will vary continuously as a function of $x$. By symmetry, $P(\mu < 0|X = 0) = P(\mu > 0|X = 0) < \frac{1}{2}$ and so there will be a range of values of $x$ containing 0 for which the posterior median is 0. By symmetry and the monotonicity of $\hat{\mu}$ there will be some threshold $t(w)$ such that the posterior median is zero if and only if $-t \leq x \leq t$. The posterior median of $\mu$, given $X = \tau$, is necessarily zero, so $\tau(w) < t(w)$.

Suppose $x > 0$. By the assumption that $\gamma$ is symmetric and unimodal, $\gamma(x - v) \geq \gamma(x + v)$ for all $v \geq 0$. Hence, multiplying by $\varphi(v)/g(x)$, if $x > 0$,

(36)   $$f(x - v|X = x, \mu \neq 0) \geq f(x + v|X = x, \mu \neq 0) \qquad \text{for all } v \geq 0.$$

Integrating over $0 < v < \infty$,

$$P(\mu \leq x|X = x, \mu \neq 0) \geq P(\mu > x|X = x, \mu \neq 0).$$

Therefore,

$$P(\mu > x|X = x) \leq P(\mu > x|X = x, \mu \neq 0) \leq \tfrac{1}{2},$$

and so the posterior median satisfies $\hat{\mu}(x) \leq x$ for all $x > 0$. By the monotonicity of $\hat{\mu}$ we have the shrinkage property $0 \leq \hat{\mu}(x) \leq x$ for all $x > 0$; by symmetry it is also the case that $0 \geq \hat{\mu}(x) \geq x$ for $x < 0$.

Finally we show that the maximum amount of shrinkage is appropriately bounded: the approach is to find a constant $a$ such that for all sufficiently large $x$,

$$(37) \qquad \begin{aligned} P(\mu > x - a | X = x) \\ = P(\mu > x - a | X = x, \mu \neq 0) P(\mu \neq 0 | X = x) > \tfrac{1}{2}. \end{aligned}$$

The term $P(\mu > x - a | X = x, \mu \neq 0)$ does not depend on $w$, and we consider it first. Set $B = \sup_{|u| \leq M} \gamma(u) e^{\Lambda u} / \gamma(M) e^{\Lambda M}$. For $u \leq 0$ and for $u \geq M$, $u \to \gamma(u) e^{\Lambda u}$ is increasing and so for any $c > M$ we have

$$(38) \qquad \begin{aligned} & \text{Odds}(\mu > c | X = x, \mu \neq 0) \\ & = \frac{\int_c^\infty \gamma(u) \varphi(x-u)\,du}{\int_{-\infty}^c \gamma(u) \varphi(x-u)\,du} \geq \frac{\int_c^\infty e^{-\Lambda u} \varphi(u-x)\,du}{B \int_{-\infty}^c e^{-\Lambda u} \varphi(u-x)\,du}. \end{aligned}$$

Since $\Lambda < \rho$, we have $\int_{-\infty}^\infty e^{-\Lambda v} \varphi(v)\,dv < \infty$, and so there is a value $a$ such that

$$(39) \qquad \int_{-a}^\infty e^{-\Lambda v} \varphi(v)\,dv > 3B \int_{-\infty}^{-a} e^{-\Lambda v} \varphi(v)\,dv.$$

As long as $x > a + M$, from (38) and (39) we will then have

$$(40) \qquad \text{Odds}(\mu > x - a | X = x, \mu \neq 0) \geq \frac{\int_{-a}^\infty e^{-\Lambda v} \varphi(v)\,dv}{B \int_{-\infty}^{-a} e^{-\Lambda v} \varphi(v)\,dv} > 3,$$

so that

$$(41) \qquad P(\mu > x - a | X = x, \mu \neq 0) > \tfrac{3}{4}.$$

Now set $\varepsilon = (\rho - \Lambda)/2$. Taking into account (26), (27) and (31), choose $\tau_1 \geq M$ large enough so that for $|u| \geq \tau_1$ we have

$$(42) \qquad (\log g)'(u) \geq -\Lambda - \varepsilon, \qquad (\log \varphi)'(u) \leq -\rho.$$

Choose $w_1$ so that $\tau(w_1) = \tau_1$, and define $c_1 = 2(\rho - \Lambda)^{-1} \log 2$. Suppose $w \leq w_1$, so that $\tau(w) \geq \tau_1$. It follows from (35) and (42) that, if $x > \tau(w) + c_1$, then

$$(43) \qquad \text{Odds}(\mu \neq 0 | X = x) = \Omega(x, w) \geq \Omega(\tau(w), w) e^{(\rho - \lambda)(x - \tau)/2} \geq 2.$$

On the other hand, if $w > w_1$ we will have $\Omega(x, w) > \Omega(x, w_1) \geq 2$ as long as $x > \tau_1 + c_1$. In either case, it follows that $P(\mu \neq 0 | X = x) \geq \tfrac{2}{3}$.

Combining this bound with (41), it follows that (37) is guaranteed whenever $x \geq \max\{a + M, \tau(w) + c_1, \tau_1 + c_1\}$; otherwise, all we can say is that $x - \hat{\mu}(x) \leq x$. Hence, for all $x > 0$ and $w \in [0, 1]$,

$$x - \hat{\mu}(x) \leq \max\{a, a + M, \tau(w) + c_1, \tau_1 + c_1\} \leq \tau(w) + b$$

with $b = \tau_1 + a \vee c_1$, which yields the required shrinkage bound since $\tau(w) \leq t(w)$.  $\square$

5.3. *Properties of posterior median for Gaussian errors.* For the remainder of the paper we specialize to the Gaussian error density $\phi$ in model (1) and to the global boundedness assumption (7) on the logarithmic derivative of the prior $\gamma$. Property (31) is then strengthened to

$$\tag{44} \sup \left| \frac{d}{du} \log g(u) \right| \leq \Lambda < \infty.$$

When $\varphi = \phi$, the representation (33) yields

$$\tag{45} 1 + \beta(y) = (g/\phi)(y) = 2 \int_0^\infty \cosh(yt) e^{-t^2/2} \gamma(t) \, dt.$$

Since $\cosh yt$ is an even convex positive function for each $t$, it follows that $1 + \beta(y)$ is also. Also, from (45) $0 < 1 + \beta(0) < 1$, so that $-1 < \beta(0) < 0$. We denote by $\beta^{-1}$ the positive inverse of $\beta$, defined on the interval $[\beta(0), \infty)$. We also have the following simple bounds:

$$\tag{46} |\beta(y)| \leq C(g/\phi)(y), \qquad \text{for all } y,$$

$$\tag{47} \tfrac{1}{2}(g/\phi)(y) \leq \beta(y) \leq (g/\phi)(y), \qquad \text{if } y > \beta^{-1}(1).$$

For a lower bound on the second derivative of $\beta$, from (45) we have, for $y > 0$,

$$\tag{48} \beta''(y) = \int_{-\infty}^\infty t^2 \cosh(yt) e^{-t^2/2} \gamma(t) \, dt > \int_{-\infty}^\infty t^2 e^{-t^2/2} \gamma(t) \, dt = \beta''(0) > 0.$$

We now develop an explicit form for the equation defining the threshold $t(w)$. Define $g_+(x) = \int_0^\infty \phi(x - \mu)\gamma(\mu) \, d\mu$ and $g_-(x) = \int_{-\infty}^0 \phi(x - \mu)\gamma(\mu) \, d\mu$. Then

$$P(\mu > 0 | X = x) = \frac{wg_+(x)}{(1 - w)\phi(x) + wg(x)}.$$

Therefore, the threshold $t$ satisfies

$$\tag{49} 2wg_+(t) = (1 - w)\phi(t) + wg(t).$$

Dividing by $w\phi(t)$ and rearranging yields

$$\tag{50} \frac{1}{w} = 1 + \frac{g_+(t) - g_-(t)}{\phi(t)} = 1 + 2 \int_0^\infty \sinh(t\mu) e^{-\mu^2/2} \gamma(\mu) \, d\mu.$$

This equation shows that the posterior median threshold $t(w)$ is continuous and strictly decreasing from $\infty$ at $w = 0$ to zero at $w = 1$.

5.4. *The link between threshold and pseudothreshold.* It will be useful to find bounds on the threshold of the posterior median function in terms of the weight $w$. It will be convenient to define the *pseudothreshold* $\zeta(w)$ by

$$\zeta = \beta^{-1}(w^{-1}).$$

The following result sets out relations between the pseudothreshold $\zeta(w)$ and the true threshold $t(w)$ of the posterior median function. In most of our discussion the dependence of $t$ and $\zeta$ on $w$ is not expressed explicitly, and, indeed, any two of $t$, $\zeta$ and $w$ can be regarded as functions of the third.

LEMMA 3.    *For all* $w \in (0,1]$,

$$(51) \qquad\qquad 1 + \beta\{t(w)\} < \beta\{\zeta(w)\} < 2 + \beta\{t(w)\}.$$

PROOF.    The bounds are a straightforward consequence of (50) defining $t(w)$, which may be rewritten in the form

$$\beta(\zeta) = \beta(t) + 2 - 2g_-(t)/\phi(t).$$

Clearly,

$$0 < g_-(t) = \int_{-\infty}^0 \phi(t-\mu)\gamma(\mu)\,d\mu < \phi(t)\int_{-\infty}^0 \gamma(\mu)\,d\mu = \tfrac{1}{2}\phi(t).$$

Thus, $0 < 2g_-(t)/\phi(t) < 1$, which establishes (51).   □

From the properties of $\beta$ we can derive two important corollaries. First, we have $0 \le t < \zeta$ for all finite $t$ and $\zeta$, so that

$$(52) \qquad\qquad\qquad t^2 < \zeta^2.$$

Second, from the property (48) that $\beta''(y) > C$ for all $y$, it follows that, for $y > 0$, $\beta'(y) > Cy$. Therefore,

$$\frac{d}{du}\beta(\sqrt{u}) = \frac{1}{2}u^{-1/2}\beta'(\sqrt{u}) > \frac{1}{2}C,$$

so that

$$\beta(\zeta) - \beta(t) = \beta(\sqrt{\zeta^2}) - \beta(\sqrt{t^2}) > \tfrac{1}{2}C(\zeta^2 - t^2).$$

Therefore,

$$(53) \qquad\qquad \zeta^2 - t^2 < 2C^{-1}\{\beta(\zeta) - \beta(t)\} \le 4C^{-1},$$

so (for a different value of $C$) $-t^2 < -\zeta^2 + C$ and so finally, for some constant $C > 0$,

$$(54) \qquad\qquad\qquad \phi(t) < C\phi(\zeta).$$

5.5. *Properties of the posterior mean.* In this section we consider the effects of using the posterior mean as an estimate instead of the posterior median. We begin by considering the behavior of the posterior distribution conditional on $\mu \neq 0$, which is also the unconditional case $w = 1$. Given any $x$, define

$$\tilde{\mu}_1(x) = E_{\text{post}}(\mu | X = x, \mu \neq 0) = \frac{\int_{-\infty}^{\infty} u\phi(x - u)\gamma(u)\,du}{\int_{-\infty}^{\infty} \phi(x - u)\gamma(u)\,du}.$$

A simple argument using the property $\phi'(t) = -t\phi(t)$ shows that $\tilde{\mu}_1(x) = x + (\log g)'(x)$, and, hence, using the bound (44),

(55)                                $|\tilde{\mu}_1(x) - x| \leq \Lambda.$

Defining $\tilde{\mu}(x, w)$ to be the posterior mean $E(\mu | X = x)$, we then have

(56)        $\tilde{\mu}(x, w) = P(\mu \neq 0 | X = x)E(\mu | X = x, \mu \neq 0) = \tilde{w}(x, w)\tilde{\mu}_1(x).$

From (36), if $x > 0$, the posterior mean $\tilde{\mu}_1(x) \leq x$; by a similar argument, for $v > 0$,

$$f_\mu(v | X = x, \mu \neq 0) > f_\mu(-v | X = x, \mu \neq 0)$$

and so $\tilde{\mu}_1(x) > 0$. Also, by a simple extension of the corresponding argument at the beginning of Section 5.2, $\tilde{\mu}_1$ is an increasing function of $x$. Hence, $\tilde{\mu}_1$ is a shrinkage rule, and from (56), so is $\tilde{\mu}(\cdot, w)$.

For each $x$ the the posterior weight $\tilde{w}(x, w)$ is monotone increasing in $w$; for $x > 0$ it follows from (56) that so also is the posterior mean $\tilde{\mu}(x, w)$.

**Bounded shrinkage properties of the posterior mean.** From (55) and (56) we have

$$x - \tilde{\mu}(x, w) = (1 - \tilde{w})x - \tilde{w}(\log g)'(x).$$

Choose $w_1$ so that $\tau(w_1) = \Lambda$, and let $\tau_2 = \tau(w \wedge w_1) \geq \Lambda$. Using (35) and (44), we have

$$\Omega(x) \geq \exp\left\{\int_{\tau_2}^x (u - \Lambda)\,du\right\} \geq \exp\left\{\int_{\tau_2}^x (u - \tau_2)\,du\right\} = \exp\left\{\tfrac{1}{2}(u - \tau_2)^2\right\}.$$

From (34), $\Omega(\tau_2) = \Omega(\tau(w \wedge w_1), w) \geq 1$ and so for $x > \tau_2$,

$$1 - \tilde{w} \leq 1/\Omega(x) \leq \exp\{-\tfrac{1}{2}(x - \tau_2)^2\}.$$

Combining this with bound (44), we obtain for $x > \tau_2$,

$$\begin{aligned}
x - \tilde{\mu}(x, w) &\leq (x - \tau_2)(1 - \tilde{w}) + \tau_2 + |(\log g)'(x)| \\
&\leq (x - \tau_2)\exp\{-\tfrac{1}{2}(x - \tau_2)^2\} + \tau_2 + \Lambda \leq e^{-1/2} + 2\Lambda + t(w).
\end{aligned}$$

If $0 \leq x \leq \tau_2$, then trivially $x - \tilde{\mu} \leq x < \tau_2 < \Lambda + t(w)$, so that we have shown that the posterior mean is a bounded shrinkage rule relative to the threshold $t(w)$.

5.6. *Bounds for integrals of exponential growth.*

LEMMA 4.   *If* $(\log h)'(z) \geq (\log k)'(z)$ *for* $z \in [\zeta_1, \zeta]$, *then*

$$\{h(\zeta)\}^{-1} \int_{\zeta_1}^{\zeta} h(z)\, dz \leq \{k(\zeta)\}^{-1} \int_{\zeta_1}^{\zeta} k(z)\, dz$$

$$\leq \begin{cases} \gamma^{-1}[1 - e^{-\gamma(\zeta - \zeta_1)}], & \text{if } k(z) = e^{\gamma z}, \\ 4(\gamma\zeta)^{-1}, & \text{if } k(z) = e^{\gamma z^2 - \alpha z}, \end{cases}$$

*where in the second case we require also that* $\zeta_1 \geq 0$ *and* $\gamma\zeta \geq \max(\alpha, 0)$.

PROOF.   The first inequality is seen easily by writing $h(z)/h(\zeta) = \exp\{-\int_z^{\zeta} (\log h)'\}$, applying the assumed inequality and integrating. The second inquality for $k(z) = e^{\gamma z}$ is trivial. For $k(z) = e^{\gamma z^2 - \alpha z}$, we first note that change of scale shows that it suffices to prove the bound for $\gamma = 1/2$. Replacing $\zeta_1$ by $0$ and completing the square, we find that the desired bound is implied by

$$(57) \qquad\qquad \int_{-\alpha}^{\zeta - \alpha} e^{v^2/2}\, dv \leq \frac{4}{\zeta} e^{(\zeta - \alpha)^2/2}.$$

If $\zeta > 2\max(\alpha, 0)$, then $\alpha < \zeta - \alpha$ and the integral on the left is bounded by $2\int_0^{\zeta - \alpha} e^{v^2/2}\, dv$. Equation (57) now follows from the inequalities

$$e^{-w^2/2} \int_0^w e^{v^2/2}\, dv = \int_0^w e^{-(w-v)(w+v)/2}\, dv$$

$$\leq \int_0^w e^{-(w-v)w/2}\, dv \leq \int_0^\infty e^{-xw/2}\, dx = 2/w. \qquad \square$$

COROLLARY 1.   *If* $g = \gamma \star \phi$ *and* $\gamma$ *satisfies* (7), *then*

$$(58) \qquad \int_0^\zeta (g/\phi)^q(x)\phi(x - \mu)\, dx \leq H_q(\zeta; \Lambda, \mu)(g/\phi)^q(\zeta)\phi(\zeta - \mu),$$

*where*

$$(59) \qquad H_q(\zeta; \Lambda, \mu) = \begin{cases} 8/[(q-1)\zeta], & \text{if } q > 1,\ \zeta > 2q\Lambda/(q-1), \\ \zeta, & \text{if } q = 1,\ \mu > \Lambda, \\ (e^{\Lambda\zeta} - 1)/\Lambda, & \text{if } q = 1,\ 0 \leq \mu \leq \Lambda. \end{cases}$$

PROOF.   Let $h(x) = (g/\phi)^q(x)\phi(x - \mu)$. Then

$$(\log h)'(x) = q(\log g)'(x) + qx - (x - \mu) \geq \begin{cases} \mu - \Lambda, & q = 1, \\ (q-1)x + \mu - q\Lambda, & q > 1. \end{cases}$$

If $q = 1$, we apply the preceding lemma with $\log k(z) = (\mu - \Lambda)z$ and $\zeta_1 = 0$ and obtain factor $H_1(\zeta; \Lambda, \mu)$ according as $\mu > \Lambda$ or not. For $q > 1$, we use the version with $\log k$ quadratic, $\gamma = (q-1)/2$ and $\alpha = q\Lambda - \mu$, so that $\gamma\zeta \geq \max(\alpha, 0)$ becomes $\zeta \geq (2/(q-1)\max(q\Lambda - \mu, 0)$.   $\square$

**6. Risk properties of thresholding procedures.** In this section we study the risk behavior of thresholding procedures. Because the thresholds obtained by the empirical Bayes procedure are data-dependent, some care is appropriate in deriving the risk. We begin with risk bounds for hard thresholding using fixed, nonrandom thresholds. These lead to comparison inequalities and so to bounds for the risk for general random thresholds. The latter continue to hold if the threshold is replaced by a pseudothreshold that is easier to find for the mixture prior model. Analogs for the posterior mean are studied as well.

6.1. *Risk bounds for fixed thresholds.* As a tool for later work, we develop risk bounds for hard thresholding,

$$\hat{\mu}_{\mathrm{HT}}(x, t) = x I\{|x| \geq t\}$$

in $L_q$ error for $0 < q \leq 2$. For the posterior mean estimator $\tilde{\mu}(x, w)$ of (56), a bound of similar structure holds for $q > 1$, based on the pseudothreshold $\zeta = \beta^{-1}(w^{-1})$ in place of $t$.

PROPOSITION 1. (a) *Fix $q \in (0, 2]$. There exists a constant $c_q \leq 4$ such that for $t \geq \sqrt{2}$ and for all $\mu$,*

$$(60) \qquad E|\hat{\mu}_{\mathrm{HT}}(X, t) - \mu|^q \leq c_q[|\mu|^q + t^{q-1}\phi(t)].$$

(b) *Now suppose $q \in (1, 2]$. There exists a constant $c'_q$ such that for $\zeta \geq \zeta(\gamma)$ and all $\mu$,*

$$(61) \qquad E|\tilde{\mu}(X, w) - \mu|^q \leq c'_q[|\mu|^q + \zeta^{q-1}\phi(\zeta)].$$

The main use of these bounds is to control risks when $\mu$ is not too large, and especially when $\mu \to 0$. The second term in each bound is, up to constants, a sharp representation of the risk at $\mu = 0$ as a function of $t$ or $\zeta$.

REMARK 2. If $q = 1$, it can be shown that the risk for the posterior mean at zero,

$$E|\tilde{\mu}(Z, w)| \geq cw \geq c\phi(\zeta)/g(\zeta),$$

is already of larger order than in (61), and so our methods for the analysis of the behavior of the posterior mean cannot immediately be extended beyond the range $1 < q \leq 2$. More remarks will be made in Section 10.

PROOF OF PROPOSITION 1. We begin with a simple bound valid for any shrinkage rule $\hat{\mu}(x)$. Indeed, for any $\mu$ and $x$,

$$(62) \qquad |\hat{\mu}(x) - \mu|^q \leq \max\{|\mu|^q, |x - \mu|^q\} \leq |\mu|^q + |x - \mu|^q.$$

Hence, if $X \sim N(\mu, 1)$ and $a_q = E|Z|^q$,

$$(63) \qquad\qquad E|\hat{\mu}(X) - \mu|^q \le |\mu|^q + a_q.$$

From this we immediately have, when $|\mu| \ge 1$,

$$(64) \qquad\qquad E|\hat{\mu}(X) - \mu|^q \le (1 + a_q)|\mu|^q \le 2|\mu|^q,$$

and so, for the rest of the proof, we confine attention to $|\mu| \le 1$, and, indeed, in view of symmetry of the risk functions, to $0 \le \mu \le 1$.

(a) Write $r_q(\mu, t)$ for the risk $E|XI\{|X| > t\} - \mu|^q$ of hard thresholding. We have

$$r_q(\mu, t) = \mu^q[\Phi(t - \mu) - \Phi(-t - \mu)] + \left(\int_{t-\mu}^{\infty} + \int_{-\infty}^{-t-\mu}\right)|z|^q\phi(z)\,dz.$$

By partial integration we obtain the upper bound

$$(65) \qquad\qquad r_q(0, t) = 2\int_t^{\infty} z^q\phi(z)\,dz \le b_q t^{q-1}\phi(t),$$

where $b_q$ may be taken as 2 for $q \le 1$ and as 4 when $1 < q \le 2$ and $t \ge \sqrt{2}$. By subtraction,

$$(66) \qquad r_q(\mu, t) - r_q(0, t) = \mu^q[\Phi(t - \mu) - \Phi(-t - \mu)] + \Delta(\mu, t),$$

where

$$\Delta(\mu) = \Delta(\mu, t) = \left(\int_{t-\mu}^{t} - \int_t^{t+\mu}\right)z^q\phi(z)\,dz.$$

The function $\phi_q(t) = t^q\phi(t)$ is positive on $(0, \infty)$ for all $q$, and for $q \ge 0$ attains its maximum value $\phi_q^* = \phi(0)(q/e)^{q/2}$ at $t = \sqrt{q}$. We remark that $\phi_q^* \le 1/2$ when $0 \le q \le 3$. Set $\phi_q(t, \mu) = \phi_q(t - \mu) + \phi_q(t + \mu)$: some calculation then shows that when $0 \le \mu \le t$,

$$\Delta''(\mu) = \phi_{q+1}(t, \mu) - q\phi_{q-1}(t, \mu) \le \phi_{q+1}(t, \mu) \le 2\phi_{q+1}^*.$$

Since $\Delta(0) = \Delta'(0) = 0$, we therefore have $\Delta(\mu) \le \phi_{q+1}^* \mu^2$, at least for $0 \le \mu \le t$. Combining this with (66), we have

$$r_q(\mu, t) \le r_q(0, t) + \mu^q + \phi_{q+1}^* \mu^2.$$

If $\mu \le 1$, then $\mu^2 \le \mu^q$, and bringing in (65), we obtain (60) with $c_q \le 4$.

(b) In view of (62) and (65), we have, for $0 \le \mu \le \frac{1}{2}$ and $\zeta \ge 2$,

$$
\begin{aligned}
E[|\tilde{\mu} - \mu|^q, |X| \ge \zeta] &\le 2\int_{\zeta}^{\infty}[\mu^q + |x - \mu|^q]\phi(x - \mu)\,dx \\
&\le 2[\mu^q(\zeta - \mu)^{-1} + 2(\zeta - \mu)^{q-1}]\phi(\zeta - \mu) \\
&\le 5\zeta^{q-1}\phi(\zeta - \mu).
\end{aligned}
$$
(67)

On the interval $0 \leq x \leq \zeta$ we have here $1/c_0 = g(0)/\phi(0) \leq 1 + w\beta(x) \leq 2$, and so

$$\tilde{w} = \frac{wg(x)/\phi(x)}{1 + w\beta(x)} \leq c_0 wg(x)/\phi(x).$$

Together with (55) this shows that $\tilde{\mu}(x, w) \leq c_0 w(\zeta + \Lambda)(g/\phi)(x)$, and so

(68)
$$E[|\tilde{\mu} - \mu|^q, |X| \leq \zeta]$$

$$\leq 2^{q-1}\mu^q + [2c_0 w(\zeta + \Lambda)]^q \int_0^\zeta (g/\phi)^q(x)\phi(x - \mu)\, dx$$

(69)
$$\leq 2\mu^q + [8/(q-1)\zeta][3c_0 w\zeta(g/\phi)(\zeta)]^q \phi(\zeta - \mu),$$

using (58) and (59), valid for $\zeta > 2q\Lambda/(q-1)$, and noting that for such $\zeta$, we have also $\zeta + \Lambda \leq 3\zeta/2$. Since $\beta(\zeta) = w^{-1} \geq 1$, we always have $(g/\phi)(\zeta) = \beta(\zeta) + 1 \leq 2\beta(\zeta) = 2w^{-1}$. Inserting these remarks into (69) and combining with (67) yields, for $\mu \in [0, \frac{1}{2}]$ and $\zeta \geq \zeta_0 = \max\{2, \beta^{-1}(1), 2q\Lambda/(q-1)\}$, that

$$E|\tilde{\mu} - \mu|^q \leq 2\mu^q + c_q\zeta^{q-1}\phi(\zeta - \mu).$$

For $0 < \mu < 1/\zeta$, one has $\phi(\zeta - \mu) \leq e\phi(\zeta)$, while for $1/\zeta < \mu < \frac{1}{2}$, some calculus shows that $\zeta^{q-1}\phi(\zeta - \mu) \leq \mu^2 \leq \mu^q$. This completes the proof of (61) for $0 \leq \mu \leq \frac{1}{2}$, while for $\mu > \frac{1}{2}$ the bound follows by a simple modification of (64). $\square$

### 6.2. *Risk bounds for general random thresholds.*

We begin with a simple bound. Suppose that $\delta$ is a shrinkage rule with the bounded shrinkage property, and that $\hat{t}$ is a random threshold with $\hat{t} \leq t$ with probability one on the event $A$. Then

(70)
$$E|\delta(X, \hat{t}) - \mu|^q I_A \leq 2E[|\delta(X, \hat{t}) - X|^q + |X - \mu|^q]I_A$$

$$\leq 2\{|t + b|^q P(A) + [E|X - \mu|^{2q}]^{1/2} P(A)^{1/2}\}$$

$$\leq 4\{t^q + b^q + 1\}P(A)^{1/2}.$$

[We have used $E|Z|^{2q} \leq (EZ^4)^{2q/4} \leq 3$ for $q \leq 2$.]

We now consider more specific risk bounds for random thresholds. The first will be particularly useful for small values of the true mean $\mu$, in conjunction with a constant $t$ which is with high probability a lower bound for the threshold $\hat{t}$.

LEMMA 5. (a) *Suppose that $0 < q \leq 2$, that $X \sim N(\mu, 1)$ and that $\hat{t}$ is a random threshold that may depend both on $X$ and on other data. Suppose that $\delta$ is a thresholding rule with the bounded shrinkage property, and let*

$$\hat{\mu} = \delta(X, \hat{t}).$$

*Suppose that $t \geq \sqrt{2}$. Then for all $\mu$,*

$$(71) \qquad E|\hat{\mu} - \mu|^q \leq c_q[|\mu|^q + t^{q-1}\phi(t) + (t^q + b^q + 1)\{P(\hat{t} < t)\}^{1/2}].$$

(b) *If $1 < q \leq 2$, a similar result holds for the posterior mean with estimated pseudothreshold $\zeta$. For $\zeta \geq \zeta_0(q, \gamma)$ and for all $\mu$,*

$$(72) \qquad E|\tilde{\mu}(x, \hat{w}) - \mu|^q \leq c'_q[|\mu|^q + \zeta^{q-1}\phi(\zeta) + (\zeta^q + b^q + 1)\{P(\hat{\zeta} < \zeta)\}^{1/2}].$$

PROOF. The method for both parts is essentially identical, so we concentrate on the thresholding case (a). Denote by $\mu^*(X)$ the effect of applying to $X$ the hard thresholding rule with threshold $t$. If $\hat{t}$ is a data dependent threshold with $\hat{t} \geq t$, then it follows from the shrinkage and thresholding properties of $\hat{\mu}$ and $\mu^*$ that both

$$(73) \qquad \text{sign}(\hat{\mu}) = \text{sign}(\mu^*) \quad \text{and} \quad 0 \leq |\hat{\mu}| \leq |\mu^*|.$$

Hence,

$$|\hat{\mu} - \mu|^q \leq \max\{|\mu|^q, |\mu^* - \mu|^q\} \leq |\mu|^q + |\mu^* - \mu|^q.$$

If we remove the overall constraint that $\hat{t} \geq t$, it remains the case that

$$|\hat{\mu} - \mu|^q I[\hat{t} \geq t] \leq |\mu|^q + |\mu^* - \mu|^q.$$

The inequality (60) for the risk of the hard thresholding rule $\mu^*$ with fixed threshold $t$ shows that for $t \geq \sqrt{2}$,

$$(74) \qquad E|\hat{\mu} - \mu|^q I[\hat{t} \geq t] \leq c_q[|\mu|^q + t^{q-1}\phi(t)].$$

Now consider the case $\hat{t} < t$. By the bounded shrinkage property and (70) it follows that

$$(75) \qquad E|\hat{\mu} - \mu|^q I[\hat{t} < t] \leq 4(t^q + b^q + 1)\{P(\hat{t} < t)\}^{1/2}.$$

Putting together the two bounds (74) and (75) completes the proof of (71).

For the posterior mean, we use the pseudothreshold $\zeta$ and set $\hat{\mu}(x) = \tilde{\mu}(x, w(\hat{\zeta}))$ and $\mu^*(x) = \tilde{\mu}(x, w(\zeta))$ in the above argument. The key monotonicity property (73) follows from that of $w \rightarrow \tilde{\mu}(x, w)$, and the analog of (74) uses (61). Finally, we use the bounded shrinkage property of the posterior mean. □

The second lemma will be used in practice for larger values of $\mu$.

LEMMA 6. *Make the same assumptions as Lemma 5 but relax the condition that $\delta$ is necessarily a strict thresholding rule; it is still required that $\delta$ has the bounded shrinkage property. Suppose that $\hat{t}$ satisfies the inequality*

$$(76) \qquad \hat{t} \leq \sqrt{d \log n} \qquad \text{with probability 1.}$$

*Let $t$ be a nonrandom threshold, possibly depending on $n$. Then*

$$(77) \qquad E|\hat{\mu} - \mu|^q \leq 8[t^q + b^q + 1 + (d \log n)^{q/2}\{P(\hat{t} > t)\}^{1/2}].$$

PROOF. To prove Lemma 6, suppose first that $\hat{t} \leq t$. From (70) with $A = \{\hat{t} \leq t\}$, we have

$$(78) \qquad E|\hat{\mu} - \mu|^q I[\hat{t} \leq t] \leq 4(t^q + b^q + 1).$$

Now use (70) again, now with $A = \{\hat{t} > t\}$ and note that $\hat{t} \leq \sqrt{d \log n}$ w.p.1 on $A$, so that

$$(79) \qquad E|\hat{\mu} - \mu|^q I[\hat{t} > t] \leq 4((d \log n)^{q/2} + b^q + 1)\{P(\hat{t} > t)\}^{1/2}.$$

Combining the two results (78) and (79) completes the proof of Lemma 6. $\square$

REMARK 3. Lemma 6 applies in particular to the posterior mean rule $\hat{\mu} = \hat{\mu}(x, w(\hat{\zeta}))$ with estimated pseudothreshold $\hat{\zeta}$.

It also follows from (52) and (54) that the bounds in Lemmas 5 and 6 remain valid if thresholds $t$ are replaced by pseudothresholds $\zeta$ throughout.

**7. Moments of the score function.** In this section we derive properties of the score function $S(w)$ that will facilitate our detailed consideration of the behavior of $\hat{w}$. Suppose that $Z \sim N(0, 1)$ and define $m_1(\mu, w) = E \beta(Z + \mu, w)$ and $m_2(\mu, w) = E \beta(Z + \mu, w)^2$. We first note that

$$\frac{\partial}{\partial \mu} m_k(\mu, w) = \int_0^\infty k \beta^{k-1}(x) \beta'(x)[1 + w\beta(x)]^{-k-1}[\phi(x - \mu) - \phi(x + \mu)] \, dx.$$

For $k = 1$, this shows that $\mu \to m_1(\mu, w)$ is increasing for $\mu \geq 0$.

7.1. *The moments $\tilde{m}(w)$ and $m_k(\mu, w)$ as functions of $w$.* We give a special name to the mean zero case and study it first:

$$(80) \qquad \tilde{m}(w) := -m_1(0, w) = -2 \int_0^\infty \beta(z, w)\phi(z) \, dz.$$

LEMMA 7. *The function $w \to \tilde{m}(w)$ is nonnegative and increasing in $w \in [0, 1]$ and satisfies $\tilde{m}(0) = 0$. If $\zeta = \beta^{-1}(w^{-1})$ is the pseudothreshold discussed in Section 5.4, then*

$$(81) \qquad \tilde{m}(w) \asymp \zeta^{\kappa-1} g(\zeta) \qquad \text{as } w \to 0.$$

LEMMA 8. *Fix $\mu > 0$. The function $w \to m_1(\mu, w)$ is decreasing in $w \in [0, 1]$ and satisfies $m_1(\mu, 0) > 0$. In terms of $\zeta = \beta^{-1}(w^{-1})$, for sufficiently small $w < w_0(\gamma)$ (not depending on $\mu$) we have*

$$(82) \qquad m_1(\mu, w) \geq \tfrac{1}{2}\beta(\zeta)\tilde{\Phi}(\zeta - \mu),$$

$$(83) \qquad m_2(\mu, w) \leq Cw^{-1} m_1(\mu, w), \qquad \mu \geq 1,$$

*while*

$$(84) \qquad m_1(\zeta, w) \sim \tfrac{1}{2} w^{-1} \qquad \text{as } w \to 0.$$

PROOF OF LEMMA 7.  For each $z \neq \beta^{-1}(0)$, $\beta(z, w)$ is a decreasing function of $w$ and so $\tilde{m}(w)$ is increasing. It follows that, as $w \searrow 0$,

$$\tilde{m}(w) \searrow \tilde{m}(0) = -\int_{-\infty}^{\infty} \beta(z)\phi(z)\,dz = \int_{-\infty}^{\infty} \{\phi(z) - g(z)\}\,dz = 0.$$

To study the asymptotic behavior of $\tilde{m}(w)$ as $w \to 0$, use the property $\int_{-\infty}^{\infty} \beta(y)\phi(y)\,dy = 0$ to obtain

$$\tilde{m}(w) = 2\int_0^{\infty} \frac{w\beta(z)^2}{1 + w\beta(z)}\phi(z)\,dz.$$

Define $\zeta = \beta^{-1}(w^{-1})$. On the range $z < \zeta$, we have $w\beta(z) < 1$, so that $1 + \beta(0) \leq 1 + w\beta(z) \leq 2$. On the other hand, for $z \geq \zeta$ we have $w\beta(z) < 1 + w\beta(z) < 2w\beta(z)$. It follows that

$$(85) \qquad \tilde{m}(w) \asymp \int_0^{\zeta} w\beta(z)^2\phi(z)\,dz + \int_{\zeta}^{\infty} \beta(z)\phi(z)\,dz.$$

Appealing to (46), (58) and (59), we then have for $\zeta \geq 4\Lambda$,

$$(86) \qquad \begin{aligned} \int_0^{\zeta} \beta(u)^2\phi(u)\,du &\leq C\int_0^{\zeta} g(u)^2/\phi(u)\,du \leq 8C\zeta^{-1}g(\zeta)^2/\phi(\zeta) \\ &\sim C\zeta^{-1}\beta(\zeta)g(\zeta) = Cw^{-1}\zeta^{-1}g(\zeta), \end{aligned}$$

using (47) and assuming also that $\zeta \geq \zeta_0$. Hence, the first integral in (85) is bounded by a term of order $\zeta^{-1}g(\zeta)$.

Because $\beta(u)\phi(u) \sim g(u)$ as $u \to \infty$, the second integral in (85) is asymptotic, via (30) to

$$\int_{\zeta}^{\infty} g(u)\,du \asymp \zeta^{\kappa-1}g(\zeta).$$

This term strictly dominates the bound $\zeta^{-1}g(\zeta)$ and, therefore, we can conclude that $\tilde{m}(w)$ is bounded above and below by multiples of $\zeta^{\kappa-1}g(\zeta)$.  □

PROOF OF LEMMA 8.  Note first that the expression

$$m_1(\mu, w) = \int_{-\infty}^{\infty} \beta(t, w)\phi(t - \mu)\,dt$$

shows that $m_1(\mu, w)$ increases monotonically as $w \to 0$. The limiting value is

$$\begin{aligned} m_1(\mu, 0) &= \int_{-\infty}^{\infty} \beta(t)\phi(t - \mu)\,dt = \int_{-\infty}^{\infty} \left(\frac{g(t)}{\phi(t)} - 1\right)\phi(t - \mu)\,dt \\ &= \int_{-\infty}^{\infty} \exp\left(\mu t - \frac{1}{2}\mu^2\right)g(t)\,dt - 1 = e^{-\mu^2/2}M_g(\mu) - 1, \end{aligned}$$

where $M_g$ denotes the moment generating function of $g$. Since $g$ is the convolution of $\gamma$ and $\phi$, and $\gamma$ is symmetric,

$$e^{-\mu^2/2} M_g(\mu) = e^{-\mu^2/2} M_\gamma(\mu) M_\phi(\mu) = M_\gamma(\mu)$$

$$= 2 \int_0^\infty \cosh(\mu t) \gamma(t) \, dt > 1,$$

so that $m_1(\mu, 0) > 0$. [If $g$ has sufficiently heavy tails, then $m_1(\mu, 0)$ may be infinite.]

For sufficiently small $w$, we have

$$\int_{-\infty}^\zeta \beta(t, w) \phi(t - \mu) \, dt \geq 0, \tag{87}$$

since the limiting value of this expression is $m_1(\mu, 0) > 0$. It follows that

$$m_1(\mu, w) \geq \int_\zeta^\infty \frac{\beta(t)}{1 + w\beta(t)} \phi(t - \mu) \, dt \geq \frac{1}{2} w^{-1} \int_\zeta^\infty \phi(t - \mu) \, dt.$$

We turn to the bound on $m_2(\mu, w)$. Notice first that

$$|\beta(x, w)| \leq \begin{cases} C = |\beta(0)| / \{1 + \beta(0)\}, & \text{if } \beta(x) < 0, \\ w^{-1}, & \text{if } \beta(x) \geq 0. \end{cases} \tag{88}$$

Hence,

$$E|\beta(\mu + Z, w)| = m_1(\mu, w) + E\{|\beta(\mu + Z, w)| - \beta(\mu + Z, w)\}$$

$$\leq m_1(\mu, w) + 2C \leq C m_1(\mu, w)$$

for sufficiently small $w$ and $\mu \geq 1$, since we then have $m_1(\mu, w) \geq m_1(1, w) \geq C$. It also follows from (88) that, again for sufficiently small $w$, $|\beta(t, w)| \leq w^{-1}$ for all $t$, and so

$$m_2(\mu, w) \leq E\{\beta(\mu + Z, w)^2\} \leq w^{-1} E|\beta(\mu + Z, w)| \leq C w^{-1} m_1(\mu, w).$$

Turning finally to the proof of (84), we have

$$m_1(\zeta, w) = \int \frac{\beta(z + \zeta)}{1 + w\beta(z + \zeta)} \phi(z) \, dz = w^{-1} \int r(\zeta, z) \phi(z) \, dz,$$

where

$$r(\zeta, z) = \frac{\beta(\zeta + z)}{\beta(\zeta) + \beta(\zeta + z)} \to I\{z > 0\}$$

as $\zeta \to \infty$, since letting $O_1(\Lambda z)$ denote a quantity bounded in absolute value by $|\Lambda z|$,

$$\frac{\beta(\zeta)}{\beta(\zeta + z)} \sim \frac{g(\zeta)}{g(\zeta + z)} \frac{\phi(\zeta + z)}{\phi(\zeta)}$$

$$= \exp\{O_1(\Lambda z) - \zeta z - z^2/2\} \to \begin{cases} 0, & z > 0, \\ \infty, & z < 0. \end{cases}$$

The conclusion (84) follows from the dominated convergence theorem, since $|r(\zeta, z)| < 1$ [at least for $\zeta$ large enough that $\beta(\zeta) > 2|\beta(0)|$].  □

7.2. *The moments $m_k(\mu, w)$ as functions of $\mu$.* We shall need a series of bounds for $m_k(\mu, w)$, each successively more refined as $\mu$ is constrained to be closer to zero.

LEMMA 9. *There are constants $C_i$ such that for all $w$, defining $c$ as in* (88),

$$(89) \qquad m_1(\mu, w) \leq \begin{cases} -\tilde{m}(w) + C_1\zeta(w)\mu^2, & \text{for } |\mu| < 1/\zeta(w), \\ C_2\phi(\zeta/2)w^{-1}, & \text{for } |\mu| \leq \zeta(w)/2, \\ (w \wedge c)^{-1}, & \text{for all } \mu \end{cases}$$

*and*

$$(90) \qquad m_2(\mu, w) \leq \begin{cases} C_3\zeta(w)^{-\kappa}w^{-1}\tilde{m}(w), & \text{for } |\mu| < 1/\zeta(w), \\ C_4\zeta^{-1}\phi(\zeta/2)w^{-2}, & \text{for } |\mu| \leq \zeta(w)/2, \\ (w \wedge c)^{-2}, & \text{for all } \mu. \end{cases}$$

PROOF. We first remark that the global bounds $m_k(\mu, w) \leq (w \wedge c)^{-k}$ follow trivially from (88). We derive a bound on the behavior of $m_1(\mu, w) - m_1(0, w)$ for small $\mu \neq 0$. Assume that $|\mu| \leq \zeta^{-1}$ and that $\zeta > 2$. Then for all $y \in [-\zeta, \zeta]$,

$$(91) \qquad \phi(y - \mu) = \phi(y)\exp\left(\mu y - \tfrac{1}{2}\mu^2\right) \leq e\phi(y)$$

and

$$(92) \qquad |\phi''(y - \mu)| = |(y - \mu)^2 - 1|\phi(y - \mu) \leq c(1 + y^2)\phi(y),$$

where the absolute constant $c' < 1.25e$. Using the property that $|\beta(z, w)|$ is bounded above by $\min\{w^{-1}, \beta(z)\}$ if $\beta(z) > 0$ and by $\{1 + \beta(0)\}^{-1}|\beta(z)|$ if $\beta(z) < 0$, it follows that

$$(93) \qquad \frac{\partial^2 m_1(\mu, w)}{\partial \mu^2} \leq \int_{-\infty}^{\infty} |\beta(z, w)\phi''(z - \mu)|\, dz$$

$$(94) \qquad \leq C\int_{-\zeta}^{\zeta} |\beta(z)|(1 + z^2)\phi(z)\, dz + 2w^{-1}\int_{|z| > \zeta} \phi''(z - \mu)\, dz.$$

Since $g(z) \leq C(1 + z^2)^{-1}$ for all $z$, it follows that

$$|\beta(z)|\phi(z) \leq |g(z) - \phi(z)| \leq C(1 + z^2)^{-1}$$

and hence that

$$(95) \qquad \int_{-\zeta}^{\zeta} |\beta(z)|(1 + z^2)\phi(z)\, dz \leq C\zeta.$$

For $|\mu| \leq \zeta^{-1}$ and $\zeta > 2$ we have

$$w^{-1} \int_{|z| > \zeta} \phi''(z - \mu) \, dz$$

$$\leq 2w^{-1} \int_{\zeta}^{\infty} \phi''(z - |\mu|) \, dz$$

(96)
$$= -2w^{-1} \phi'(\zeta - |\mu|) = 2w^{-1}(\zeta - |\mu|)\phi(\zeta - |\mu|)$$

$$\leq C\zeta\beta(\zeta)\phi(\zeta) \leq C\zeta(1 + \zeta^2)^{-1}.$$

Combining (95) and (96), recalling that $C$ can be a different constant in different expressions, we can conclude that, for $|\mu| \leq \zeta^{-1}$ and $\zeta > 2$,

$$\frac{\partial^2 m_1(\mu, w)}{\partial \mu^2} \leq C\zeta.$$

Since by symmetry $\partial m_1(\mu, w)/\partial \mu = 0$ when $\mu = 0$, it follows that, for $\mu \leq \zeta^{-1}$,

$$m_1(\mu, w) - m_1(0, w) \leq C\zeta\mu^2,$$

which completes the proof of (89).

Turn now to the second moment. Suppose throughout that $|\mu| \leq \zeta^{-1}$ and $\zeta > 2$ and, without loss of generality, that $\mu \geq 0$. By the bounds on $|\beta(z)/\{1 + w\beta(z)\}|$ and on $\phi(z - \mu)$ as above,

$$m_2(\mu, w) \leq C \int_{-\zeta}^{\zeta} \beta(z)^2 \phi(z - \mu) \, dz + \int_{|z| > \zeta} \beta(\zeta)^2 \phi(z - \mu) \, dz$$

(97)
$$\leq C \int_0^{\zeta} \beta(z)^2 \phi(z) \, dz + 2\beta(\zeta)^2 \tilde{\Phi}(\zeta - \mu)$$

(98)
$$\leq Cw^{-1}\zeta^{-1} g(\zeta) + C\beta(\zeta)^2 \phi(\zeta - \mu)/(\zeta - \mu)$$

by using (86). To deal with the second term in (98), use the bounds on $\phi(\zeta - \mu)$ and the property that $\zeta - \mu > \zeta/2$ to conclude that

$$\beta(\zeta)^2 \phi(\zeta - \mu)/(\zeta - \mu) \leq C\zeta^{-1}\beta(\zeta)^2 \phi(\zeta) \leq C\zeta^{-1}\beta(\zeta)g(\zeta).$$

It follows that, for $|\mu| \leq \zeta^{-1}$ and $\zeta > 2$,

$$m_2(\mu, w) \leq C\zeta^{-1}\beta(\zeta)g(\zeta).$$

Now use the property (81) that $g(\zeta) \leq C\zeta^{1-\kappa}\tilde{m}(w)$ to complete the proof of (90).

We now turn to the proof of the intermediate bounds. Note first that

$$m_k(\mu, w) \leq 2 \int_0^{\infty} \left[ \frac{\beta(x)}{1 + w\beta(x)} \right]^k \phi(x - \mu) \, dx.$$

On $[0, \zeta]$ we have $1 + w\beta(x) \geq 1 + \beta(0) > 0$, so that $[1 + w\beta]^{-1} \leq C$, while on $[\zeta, \infty]$ clearly $w\beta/(1 + w\beta) \leq 1$. Hence

$$m_k(\mu, w) \leq C \int_0^\zeta \beta^k(x)\phi(x - \mu)\, dx + 2w^{-k} \int_\zeta^\infty \phi(x - \mu)\, dx$$

$$= CI_{k,\zeta} + 2w^{-k}I'_{k,\zeta}.$$

Since $\beta(\zeta) = w^{-1}$, we have for $|\mu| \leq (1 - a)\zeta$,

$$I'_{k,\zeta} = \tilde{\Phi}(\zeta - \mu) \leq \phi(a\zeta)/a\zeta.$$

Turning now to $I_{k,\zeta} \leq C \int_0^\zeta (g/\phi)^k(x)\phi(x - \mu)\, dx$, we apply (58) and (59): since $(g/\phi)(\zeta) \leq 2w^{-1}$ and $\phi(\zeta - \mu) \leq \exp\{-\zeta(\zeta - 2\mu)/4\}\phi(\zeta/2)$ for $0 \leq \mu \leq \zeta/2$, we have

$$I_{k,\zeta} \leq 2^k w^{-k}\phi(\zeta/2)H_k(\zeta; \Lambda, \mu)\exp\{-\zeta(\zeta - 2\mu)/4\}.$$

The desired conclusion for $k = 2$ follows. For $k = 1$ one may check that

$$\sup_{\zeta > 8\Lambda, 0 \leq \mu \leq \zeta/2} H_1(\zeta; \Lambda, \mu)\exp\{-\zeta(\zeta - 2\mu)/4\} \leq C(\Lambda),$$

while a direct argument shows that for $\zeta \leq 8\Lambda$, regardless of $\mu$, $I_{1,\zeta} \leq Cw^{-1} \leq C\phi(\zeta/2)w^{-1}$. $\quad \square$

**8. The marginal maximum likelihood weight and its risk properties.** The marginal maximum likelihood method yields a random weight $\hat{w}$, dependent on all the data $X_1, \ldots, X_n$, and, hence, to a random threshold and pseudothreshold. In this section we study the properties of $\hat{w}$ in order to use the risk bounds of Section 6.2 to bound the risk for the whole procedure and, hence, complete the proof of Theorem 1. The structure of the proof is essentially the same for both nearly black and $\ell_p$ sparseness classes, and to avoid unnecessary repetition of arguments, it is helpful to define

$$\tag{99} \tilde{p} = \begin{cases} p, & \text{if } p > 0, \\ 1, & \text{if } p = 0. \end{cases}$$

The bounds obtained in Theorems 1 and 2 do not change as $p$ increases above 2. Furthermore, for $p \geq 2$ we have $\ell_p[\eta] \subset \ell_2[\eta]$ and so demonstrating the bounds for $p = 2$ will imply that they hold for all larger $p$. Therefore, for the whole of the subsequent argument, we assume without loss of generality that $p \leq 2$.

The strategy is to consider separately the components of the risk for large and small $\mu_i$. We employ risk decompositions

$$\tag{100} \begin{aligned} R_q(\hat{\mu}, \mu) &= n^{-1} \sum_{|\mu_i| \leq \tau} E|\hat{\mu}_i - \mu_i|^q + n^{-1} \sum_{|\mu_i| > \tau} E|\hat{\mu}_i - \mu_i|^q \\ &= R_q(\tau) + \tilde{R}_q(\tau), \end{aligned}$$

say. For "global" risk bounds, we take $\tau = 1$, while for risk bounds over $\ell_p[\eta]$, we take $\tau$ roughly of order $(2 \log \eta^{-\widetilde{p}})^{1/2}$.

In each case properties of the estimated weight and corresponding pseudothreshold are derived; these are then substituted into the appropriate expression for the risk. We begin by the consideration of the threshold and risk for the components with small $\mu_i$.

8.1. *Small signals*: *lower bounds for thresholds.* Suppose that, for some $p$ with $0 < p \leq 2$ and for some $\eta > 0$, $\mu$ lies in an $\ell_p$ ball:

$$\ell_p[\eta] = \left\{ \mu : n^{-1} \sum_{i=1}^{n} |\mu_i|^p \leq \eta^p \right\}$$

or, if $p = 0$, that the proportion of $\mu_i$ that are nonzero is at most $\eta$. Let $Z_i$ be independent $N(0,1)$ random variables, and let $\hat{w}$ be the weight estimated from the data $\mu_i + Z_i$ by the marginal maximum likelihood procedure. Define the pseudothreshold $\hat{\zeta} = \beta^{-1}(\hat{w}^{-1})$.

One cannot hope to adapt to signals that are too small relative to the sample size $n$; this corresponds to restricting $t(\hat{w})$ to the range $[0, \sqrt{2 \log n}]$. Hence, we set

$$\tilde{\eta}^{\widetilde{p}} = \max \{\eta^{\widetilde{p}}, n^{-1} (\log n)^2\}$$

and, with the usual definition $\zeta = \beta^{-1}(w^{-1})$, define the weight $w = w(\eta, n)$ by

$$(101) \qquad \zeta^{p - \kappa} w \tilde{m}(w) = \tilde{\eta}^{\widetilde{p}}.$$

Writing the left-hand side as the product of $\zeta^{p - \kappa}/\beta(\zeta)$ and $\tilde{m}(w)$, both of which are increasing in $w$ (for $w$ sufficiently small), shows that $w$ is well defined and monotonically increasing in $\tilde{\eta}$, at least for $\tilde{\eta}$ small.

The intent of this definition is to choose a weight $w = w(\eta, n)$ and pseudo threshold $\zeta = \zeta(\eta, n)$ which is both a lower bound to $\hat{\zeta} = \zeta(\hat{w})$ for $\mu \in \ell_p[\eta]$ with high probability (Lemma 10) and is of the right size to yield minimax risk bounds [see (103), (104) and Section 8.2].

**Some properties of $w$ and $\zeta$.** Using the definition of $\beta$ and the property (81) that $\tilde{m}(w) \asymp \zeta^{\kappa - 1} g(\zeta)$,

$$(102) \qquad \tilde{\eta}^{\widetilde{p}} \asymp \zeta^{p - \kappa} w \zeta^{\kappa - 1} g(\zeta) \asymp \zeta^{p-1} [g(\zeta)/\beta(\zeta)] \asymp \zeta^{p-1} \phi(\zeta).$$

We immediately obtain a bound,

$$(103) \qquad \zeta \phi(\zeta) \asymp \tilde{\eta}^{\widetilde{p}} \zeta^{2-p}.$$

Taking logarithms,

$$|\log \tilde{\eta}^{-\widetilde{p}} - \tfrac{1}{2}\zeta^2 + (p-1)\log\zeta| < C.$$

Hence, as $\tilde{\eta} \to 0$,

(104)                                    $$\zeta^2 \sim 2\log\tilde{\eta}^{-\widetilde{p}}.$$

More explicitly, there exist constants $c$ such that

(105)   $\zeta^2 \geq \begin{cases} 2\log\eta^{-\widetilde{p}} + (p-1)\log\log\eta^{-\widetilde{p}} - c, & \text{if } \eta^{\widetilde{p}} \geq n^{-1}\log^2 n, \\ 2\log n - (5-p)\log\log n - c, & \text{if } \eta^{\widetilde{p}} \leq n^{-1}\log^2 n. \end{cases}$

Approximation (104) shows that our pseudothreshold bound $\zeta(\eta, n)$ has the order of the minimax threshold for $\ell_p[\eta_n]$, and the right-hand side of (103) is essentially the asymptotic expression for the normalized minimax risk. We now show that $\zeta(\eta, n)$ is typically a lower bound for the estimated pseudothreshold when the signal is small.

LEMMA 10.   *Let the pseudothreshold $\zeta = \zeta(\eta, n)$ corresponding to $\tilde{\eta}$ be defined by* (101). *There exist $C = C(\gamma)$ and $\eta_0 = \eta_0(\gamma)$ such that if $\eta \leq \eta_0$ and $n/\log^2 n \geq \eta_0^{-\widetilde{p}}$, then*

(106)                      $$\sup_{\mu \in \ell_p[\eta]} P_\mu(\hat{\zeta} < \zeta) \leq \exp\{-C(\log n)^{3/2}\}.$$

It follows from this lemma that if $\mu$ is very sparse ($\eta^{\widetilde{p}} \leq n^{-1}\log^2 n$), then $\hat{t}$ and $\hat{\zeta}$ are, in relative terms, close to $\sqrt{2\log n}$. On the other hand, if $\mu$ is less sparse, then $\hat{t}$ and $\hat{\zeta}$ are at least about $\sqrt{2\log\eta^{-p}}$. (Recall from (53) that the difference $\hat{\zeta} - \hat{t} \in [0, C/\hat{t})$ is small.)

PROOF OF LEMMA 10.   This argument leading to (103) also shows that

$$\tfrac{1}{2}\zeta^2 - (p-1)\log\zeta \leq \log n - 2\log\log n + O(1),$$

and hence that $t(w) < \zeta(w) \leq \sqrt{2\log n}$ for $n$ sufficiently large, so that $w \in [w_n, 1]$, the interval over which the likelihood $\ell(w)$ is maximized. Consequently, $\{\hat{\zeta} < \zeta\} = \{\hat{w} > w\} = \{S(w) > 0\}$. The summands in $S(w) = \sum_{i=1}^n \beta(\mu_i + Z_i, w)$ are independent, and in view of (88), bounded by $c_0 w^{-1}$.

We therefore recall Bernstein's inequality [e.g., Pollard (1984), page 193], which gives exponential bounds on the tail probabilities of the sum of uniformly bounded independent random variables.

PROPOSITION 2.   *Suppose that $W_1, W_2, \ldots, W_n$ are independent random variables with $EW_i = 0$ and $|W_i| \leq M$ for $i = 1, \ldots, n$. Suppose that $V \geq \sum_{i=1}^n \operatorname{var} W_i$. Then, for any $A > 0$,*

$$P\left(\sum_{i=1}^n W_i > A\right) \leq \exp\{-\tfrac{1}{2}A^2/(V + \tfrac{1}{3}MA)\}.$$

We have

$$(107) \qquad P(\hat{w} > w) = P\{S(w) > 0\} = P\left(\sum_{i=1}^{n} W_i > A\right),$$

where $W_i = \beta(\mu_i + Z_i, w) - m_1(\mu_i, w)$, $M = 2(c \wedge 1)^{-1} w^{-1}$ and

$$A = \sum_{i=1}^{n} -m_1(\mu_i, w).$$

Define sets of "small," "medium" and "large" co-ordinates,

$$\mathcal{S} = \{i : |\mu_i| \leq \zeta^{-1}\}, \qquad \mathcal{M} = \{i : \zeta^{-1} < |\mu_i| \leq \zeta/2\}, \qquad \mathcal{L} = \{i : |\mu_i| > \zeta/2\}.$$

For the nearly black case, it suffices to consider only two classes, coalescing $\mathcal{M}$ and $\mathcal{S}$, but it is notationally simpler to use the same argument as for $\ell_p$.

Using the three parts of (89),

$$\sum m_1(\mu_i, w) \leq \sum_{i \in \mathcal{S}} [-\tilde{m}(w) + C_1 \zeta \mu_i^2] + C|\mathcal{M}|\phi(\zeta/2) w^{-1} + |\mathcal{L}| w^{-1}.$$

On the $\ell_p$-ball $\ell_p[\eta]$, we have $\#\{i : |\mu_i| > t\} \leq n\eta^p t^{-p}$ and so

$$(108) \qquad |\mathcal{S}| \geq n - n\eta^{\widetilde{p}} \zeta^p, \qquad |\mathcal{M}| \leq n\eta^{\widetilde{p}} \zeta^p, \qquad |\mathcal{L}| \leq n\eta^{\widetilde{p}} 2^p \zeta^{-p}.$$

On the set $\mathcal{S}$, we have $\mu_i^2 \leq |\mu_i|^p \zeta^{p-2}$, and so, on making use of (101),

$$\begin{aligned} \sum m_1(\mu_i, w) &\leq -n\tilde{m}(w) + Cn\eta^{\widetilde{p}} \zeta^{-p} w^{-1} \\ &\qquad \times [w\zeta^{2p}\tilde{m}(w) + C_1 w\zeta^{2p-1} + \zeta^{2p}\phi(\zeta/2) + 1] \\ &\leq -n[\tilde{m}(w) - C\tilde{\eta}^{\widetilde{p}} \zeta^{-p} w^{-1}] \\ &= -n\tilde{m}(w)[1 - C\zeta^{-\kappa}] \leq -\tfrac{1}{2} n\tilde{m}(w) \end{aligned}$$

for $w < w_0$. Consequently, $A \geq \tfrac{1}{2} n\tilde{m}(w)$.

We now obtain a bound on $V = \sum \operatorname{var} W_i$. Using the same decomposition into small, medium and large coordinates, we have from the three parts of (90),

$$V \leq \sum m_2(\mu_i, w) \leq C_3 |\mathcal{S}| \zeta^{-\kappa} w^{-1} \tilde{m}(w) + C|\mathcal{M}| w^{-2} \zeta^{-1} \phi(\zeta/2) + |\mathcal{L}| w^{-2}.$$

Using now (108) along with (101), we find that for sufficiently small $w$,

$$\begin{aligned} V &\leq Cn\zeta^{-\kappa} w^{-1} \tilde{m}(w) + Cn\eta^{\widetilde{p}} \zeta^{p-1} w^{-2} \phi(\zeta/2) + Cn\eta^{\widetilde{p}} \zeta^{-p} w^{-2} \\ &\leq Cnw^{-1} \tilde{m}(w)[\zeta^{-\kappa} + \zeta^{2p-\kappa-1}\phi(\zeta/2) + \zeta^{-\kappa}] \\ &\leq Cnw^{-1} \tilde{m}(w) \zeta^{-\kappa}. \end{aligned}$$

Turning to the exponent in the Bernstein bound, we have for $w \leq w_0$,

$$\left[\frac{A^2}{V + (1/3)MA}\right]^{-1} = \frac{V}{A^2} + \frac{M}{3A} \leq \frac{Cnw^{-1}\tilde{m}(w)\zeta^{-\kappa}}{n^2\tilde{m}(w)^2} + \frac{Cw^{-1}}{n\tilde{m}(w)}$$
$$= C\{nw\tilde{m}(w)\}^{-1}.$$

Therefore, applying the definition of $\tilde{\eta}$,

$$\frac{A^2}{V + (1/3)MA} \geq Cnw\tilde{m}(w) \geq Cn\tilde{\eta}^{\tilde{p}}\zeta^{\kappa-p} \geq C(\log n)^2[\log\tilde{\eta}^{-\tilde{p}}]^{(\kappa-p)/2}.$$

Define $n_0(\gamma)$ so that $n \geq n_0$ if and only if $n/\log^2 n \geq \eta_0^{-\tilde{p}}$. If $\kappa > p$, then for $\eta < \eta_0$ and $n \geq n_0$ we have $\tilde{\eta}^{-\tilde{p}} = \min\{\eta^{-\tilde{p}}, n/\log^2 n\} \geq \eta_0^{-p}$, while if $\kappa < p$, then $\tilde{\eta}^{-p} \leq n$ and $\frac{\kappa-p}{2} > -1/2$. In either case we have

$$[\log\tilde{\eta}^{-\tilde{p}}]^{(\kappa-p)/2} \geq C(\log n)^{-1/2}.$$

Applying the Bernstein inequality to (107) concludes the proof of (106), so long as $w \leq w_0(\gamma)$. Use (101) to define $\tilde{\eta}_0$ as the value of $\tilde{\eta}$ corresponding to $w_0$, and then set $\eta_0 = \tilde{\eta}_0$ to arrive at the first statement of Lemma 10. $\square$

8.2. *Small signals: risk behavior.* We apply risk bound (71) of Lemma 5. Bound (54) permits the inequality to be rewritten in terms of the pseudothreshold $\zeta$. We have, for all values of $\mu_i$,

$$(109) \qquad E|\hat{\mu}_i - \mu_i|^q \leq C\{|\mu_i|^q + \zeta^{q-1}\phi(\zeta) + (1 + \zeta^q)P(\hat{\zeta} < \zeta)^{1/2}\}.$$

If $\tilde{\eta}$ is sufficiently small [less than $\eta_0 = \eta_0(\gamma)$, say], then we may use bounds (104) and (103) along with the probability bound (106) to yield

$$E|\hat{\mu}_i - \mu_i|^q \leq C\{|\mu_i|^q + \tilde{\eta}^{\tilde{p}}(\log\tilde{\eta}^{-\tilde{p}})^{(q-p)/2} + (\log\tilde{\eta}^{-\tilde{p}})^{q/2}e^{-C\log^{3/2}n}\}.$$

If $n > n_0(\gamma)$, we have $\exp\{-C(\log n)^{3/2}\} \leq n^{-1}(\log n)^{2-p/2} \leq \tilde{\eta}^{\tilde{p}}(\log\tilde{\eta}^{-\tilde{p}})^{-p/2}$, and we finally obtain

$$E|\hat{\mu}_i - \mu_i|^q \leq C\{|\mu_i|^q + \tilde{\eta}^{\tilde{p}}(\log\tilde{\eta}^{-\tilde{p}})^{(q-p)/2}\}.$$

If $\tilde{\eta}^{\tilde{p}} = n^{-1}(\log n)^2$, then $\log\tilde{\eta}^{-\tilde{p}} = \log n - 2\log\log n \sim \log n$ so that, in general,

$$\tilde{\eta}^{\tilde{p}}(\log\tilde{\eta}^{-\tilde{p}})^{(q-p)/2} \leq \max\{\eta^{\tilde{p}}(\log\eta^{-\tilde{p}})^{(q-p)/2}, Cn^{-1}(\log n)^{2+(q-p)/2}\}.$$

Combining the last two expressions and summing over $i$ yields

$$R_q(\zeta) = n^{-1}\sum_{|\mu_i|\leq\zeta}E|\hat{\mu}_i - \mu_i|^q$$

$$\leq C\left\{n^{-1}\sum_{|\mu_i|\leq\zeta}|\mu_i|^q + \eta^{\tilde{p}}(\log\eta^{-\tilde{p}})^{(q-p)/2} + n^{-1}(\log n)^{2+(q-p)/2}\right\}.$$

If $q \leq p$, we have by Hölder's inequality

$$(110) \qquad n^{-1} \sum |\mu_i|^q \leq \left( n^{-1} \sum |\mu_i|^p \right)^{q/p} \leq \eta^q,$$

and also $\eta^p (\log \eta^{-p})^{(q-p)/2} \leq c_{p,q} \eta^q$ for $\eta < e^{-1}$. If $q > p$, we have $|\mu_i|^q \leq |\mu_i|^p \zeta^{q-p}$, and so, using the property that if $p = 0$ at most $n\eta$ of the terms will be nonzero,

$$n^{-1} \sum_{|\mu_i| \leq \zeta} |\mu_i|^q \leq \eta^{\widetilde{p}} \zeta^{q-p} \leq C \eta^{\widetilde{p}} (\log \eta^{-\widetilde{p}})^{(q-p)/2}.$$

In every case then, for $\mu \in \ell_p[\eta]$, $\eta \leq \eta_0(\gamma)$ and $n \geq n_0(\gamma)$, we have

$$(111) \qquad R_q(\zeta) \leq C\{r_{p,q}(\eta) + n^{-1}(\log n)^{2+(q-p)/2}\}.$$

Before leaving the consideration of small $\mu_i$, consider the case where there are no constraints on $\mu$ at all. The application of the elementary risk bound (63), along with $a_q \leq 1$, then yields an absolute bound on the average risk for small $\mu$:

$$(112) \qquad n^{-1} \sum_{|\mu_i| \leq 1} E|\hat{\mu}_i - \mu_i|^q \leq n^{-1} \sum_{|\mu_i| \leq 1} (1 + |\mu_i|^q) \leq 2.$$

8.3. *Large signals: upper bounds for thresholds.* Define $\tilde{\pi}(\tau; \mu) = n^{-1} \#\{i : |\mu_i| \geq \tau\}$. We will be interested in deriving upper bounds on the estimated pseudothreshold $\hat{\zeta}$ when it is known that $\tilde{\pi}(\tau; \mu) \geq \pi$ for appropriate choices of $\tau$.

Choose $w_0$ small enough so that both (81) and (83) apply. Define

$$(113) \qquad w(\tau, \pi) = \sup\{w \leq w_0 : \pi m_1(\tau, w) \geq 2\tilde{m}(w)\}.$$

Since $m_1(\tau, w)/\tilde{m}(w) \to \infty$ as $w \to 0$, certainly $w(\tau, \pi)$ is well defined. On the pseudothreshold scale, we write $\zeta_{\tau,\pi}$ or $\zeta(\tau, \pi)$ for $\beta^{-1}(1/w(\tau, \pi))$.

LEMMA 11. *There exist $C = C(\gamma)$ and $\pi_0 = \pi_0(\gamma)$ such that if $\pi < \pi_0$, then for all $\tau \geq 1$,*

$$(114) \qquad \sup_{\mu \, : \, \tilde{\pi}(\tau; \mu) \geq \pi} P_\mu(\hat{\zeta} > \zeta_{\tau,\pi}) \leq \exp\{-Cn\zeta_{\tau,\pi}^{\kappa-1} \phi(\zeta_{\tau,\pi})\}.$$

PROOF. If $n\pi$ of the $\mu_i$ for which $|\mu_i| \geq \tau$ are shrunk to $\pm\tau$, and all the other $\mu_i$ are set to zero, then the distribution of each $|\mu_i + Z_i|$ will be stochastically reduced. Since $\beta(y, w)$ is an increasing function of $|y|$ for each $w$, it follows that $S(w)$ will be stochastically reduced, and so $P(S(w) < 0)$ will be, if anything, increased. Thus, the maximum value of $P(\hat{\zeta} > \zeta)$ subject to the constraint that at least $n\pi$ of the $|\mu_i|$ exceed $\tau$ will be taken when exactly $n\pi$ of the $|\mu_i|$ are equal to $\tau$ and the remainder are zero. We shall therefore assume that this is the case.

We now return to the problem of bounding the probability that $S(w)$ is negative, for $w = w(\tau, \pi)$. We have, following (107) but changing the sign,

$$P(\hat{w} < w) = P(S(w) < 0) = P\left\{ \sum_{i=1}^n W_i > A \right\},$$

where, on this occasion,

$$W_i = m_1(\mu_i, w) - \beta(\mu_i + Z_i, w) \quad \text{and} \quad A = \sum_{i=1}^n m_1(\mu_i, w).$$

Just as above, $|W_i| \le 2c_0 w^{-1}$ for all $i$. To obtain a bound on $A$, we have, making use of the definition (113) of $w$,

$$n^{-1}A = (1 - \pi)m_1(0, w) + \pi m_1(\tau, w)$$
$$\ge -\tfrac{1}{2}\pi m_1(\tau, w) + \pi m_1(\tau, w) = \tfrac{1}{2}\pi m_1(\tau, w).$$

We now seek an upper bound on the sum of the variances of the $W_i$. Making use of the bound (90) for $m_2(0, w)$, bound (83) for $m_2(\tau, w)$ and (113),

$$n^{-1}\sum_{i=1}^n \operatorname{var} W_i \le m_2(0, w) + \pi m_2(\tau, w)$$
$$\le C\zeta(w)^{-\kappa} w^{-1}\tilde{m}(w) + Cw^{-1}\pi m_1(\tau, w)$$
$$\le Cw^{-1}\pi m_1(\tau, w).$$

Substituting into the expression needed for the application of Bernstein's inequality, we have

$$n\left( \frac{V}{A^2} + \frac{M}{3A} \right) \le Cw^{-1}\pi^{-1}m_1^{-1}(\tau, w),$$

so that

$$\frac{A^2}{V + (1/3)MA} \ge Cnw\pi m_1(\tau, w)$$
$$\ge Cnw\tilde{m}(w) \ge Cn\zeta^{\kappa-1}\beta(\zeta)^{-1}g(\zeta)$$
$$\ge Cn\zeta^{\kappa-1}\phi(\zeta),$$

(since $w \le w_0$).   $\square$

8.4. *Large signals*: *risk behavior*. Let $\zeta = \zeta(\tau, \tilde{\pi}(\tau; \mu))$, where $\tau$ remains unspecified for the moment. For each $\mu_i$, we have from (77) and (52),

$$E|\hat{\mu}_i - \mu_i|^q \leq C\{1 + \zeta^q + (\log n)^{q/2} P(\hat{\zeta} \geq \zeta)^{1/2}\}.$$

We then consider two cases. If $\zeta^2 > \log n$, then the right-hand side is bounded by $C(1 + 2\zeta^q)$. On the other hand, if $\zeta^2 \leq \log n$, then

$$n\zeta^{\kappa-1} \exp\{-\tfrac{1}{2}\zeta^2\} \geq Cn \exp\{-\tfrac{3}{4}\zeta^2\} > Cn^{1/4},$$

so that from (114),

$$(\log n)^{q/2} P(\hat{\zeta} \geq \zeta)^{1/2} < \log n \exp(-Cn^{1/4}) \leq 1$$

if $n \geq n_0$. It follows that, for sufficiently small $\pi$ and $n > n_0$, whether or not $\zeta^2 > \log n$,

$$(115) \qquad E|\hat{\mu}_i - \mu_i|^q \leq C\{1 + \zeta^q\}.$$

Hence,

$$(116) \qquad \tilde{R}_q(\tau) = n^{-1} \sum_{|\mu_i| \geq \tau} E|\hat{\mu}_i - \mu_i|^q \leq C\tilde{\pi}(\tau; \mu)[1 + \zeta^q(\tau, \tilde{\pi}(\tau, \mu))].$$

For the global risk bound needed for Theorem 1, we set $\tau = 1$. Let $\pi = \tilde{\pi}(1; \mu)$. We seek a bound for $\zeta = \zeta(1, \pi)$. Since $\tilde{m}(w) \asymp \zeta^{\kappa-1} g(\zeta)$ by (81), it follows that for sufficiently small $\pi$ and, hence, $w$,

$$\pi^{-1} = \frac{m_1(1, w)}{2\tilde{m}(w)} \geq C\zeta^{-\kappa} e^\zeta.$$

Taking logarithms, we have

$$\log \pi^{-1} \geq c - \kappa \log \zeta + \zeta$$

and, hence, for sufficiently small $\pi$,

$$(117) \qquad \zeta^q = \zeta(1, \pi)^q \leq 2^q (\log \pi^{-1})^q.$$

In combination with (116), this yields, regardless of the value of $\pi = \tilde{\pi}(1; \mu)$,

$$(118) \qquad \tilde{R}_q(1) \leq C\pi[1 + (\log \pi^{-1})^2] \leq C.$$

Write $\zeta_1$ for the pseudothreshold $\zeta(\eta, n)$ defined by (101). Our main goal now is to establish a large signal complement to inequality (111), namely,

$$(119) \qquad \tilde{R}_q(\zeta_1) \leq C\{r_{p,q}(\eta) + n^{-1}(\log n)^{2+(q-p)/2}\}.$$

The approach will be to apply Lemma 11 with $\tau = \zeta_2 = \zeta_2(\mu)$ defined by

$$(120) \qquad \zeta_2 = \zeta(\zeta_1, \pi), \qquad \pi = \tilde{\pi}(\zeta_1; \mu).$$

Let us first verify that, as one would expect for $\mu \in \ell_p[\eta]$, $\zeta_2 > \zeta_1$. Since $\tilde{m}(w) \asymp \zeta^{\kappa-1} g(\zeta)$ by (81), and $m_1(\zeta_1, w_1) \sim (2w_1)^{-1}$ by (84), we have, using (102),

$$\frac{m_1(\zeta_1, w_1)}{\tilde{m}(w_1)} \asymp \frac{1}{2} \frac{\beta(\zeta_1)}{g(\zeta_1)} \zeta_1^{1-\kappa} \asymp \frac{\zeta_1^{1-\kappa}}{\phi(\zeta_1)} \asymp \tilde{\eta}^{-\tilde{p}} \zeta_1^{p-\kappa}.$$

For $p > 0$ we now use the bound $\sum |\mu_i|^p \leq n\eta^p$, while for $p = 0$ we simply use $\pi \leq \eta$. Both cases are encompassed by the inequality

$$\pi \leq \eta^{\tilde{p}} \zeta_1^{-p} \leq \tilde{\eta}^{\tilde{p}} \zeta_1^{-p}, \tag{121}$$

and so

$$\frac{m_1(\zeta_1, w_1)}{\tilde{m}(w_1)} \leq C \zeta_1^{-\kappa} \pi^{-1} \ll 2\pi^{-1} \qquad \text{for } \zeta_1 \text{ large,}$$

which shows that $\zeta_2 > \zeta_1$ (and, in particular, that $\zeta_2 > 1$).

In this notation the bound (116) becomes

$$\tilde{R}_q(\zeta_1) \leq C\pi(1 + \zeta_2^q) \leq C\pi\zeta_2^q.$$

Although (121) places an upper bound on $\pi$, in fact, it may be arbitrarily much smaller. The analysis to follow considers separately cases in which $\pi$ is comparable to, or much smaller than, $\tilde{\eta}^{\tilde{p}} \zeta_1^{-p}$.

Recalling the lower bound (82) that $m_1(\zeta_1, w) \geq \frac{1}{2}\beta(\zeta)\tilde{\Phi}(\zeta - \zeta_1)$, it follows that $\zeta_2 \leq \zeta_3 = \zeta(w_3)$, where $w_3$ is the solution to

$$\tilde{\Phi}(\zeta(w) - \zeta_1) = 4\pi^{-1} w\tilde{m}(w). \tag{122}$$

$\zeta_3$ is intended as a more manageable version of $\zeta_2$.

Suppose first that $\zeta_3 > \zeta_1 + 1$. Then from (122),

$$\pi\zeta_2^q \leq C\zeta_3^{q+\kappa-1} \frac{g(\zeta_3)}{\beta(\zeta_3)} \frac{\zeta_3 - \zeta_1}{\phi(\zeta_3 - \zeta_1)}$$

$$\leq C\zeta_3^{\kappa+q} \frac{\phi(\zeta_3)}{\phi(\zeta_3 - \zeta_1)} = C\zeta_3^{\kappa+q} e^{-(\zeta_3 - \zeta_1)\zeta_1} \phi(\zeta_1).$$

Using (103) and the fact that $\zeta_3 \to \zeta_3^{\kappa+q} e^{-(\zeta_3 - \zeta_1)\zeta_1}$ is decreasing, at least for $\zeta_3 \geq \zeta_1 + 1$, we get

$$\pi\zeta_2^q \leq C\tilde{\eta}^{\tilde{p}} \zeta_1^{1-p}(\zeta_1 + 1)^{\kappa+q} e^{-\zeta_1}.$$

From (104), we then conclude that for $\zeta_1$ sufficiently large,

$$\pi\zeta_2^q \leq \max\{C\eta^{\tilde{p}}(2\log \eta^{-\tilde{p}})^{-2}\zeta_1^{\kappa+3+q-p} e^{-\zeta_1}, Cn^{-1}\zeta_1^{\kappa+7+q-p} e^{-\zeta_1}\}$$

$$\leq C\max\{\eta^{\tilde{p}}(2\log \eta^{-\tilde{p}})^{-2}, n^{-1}\}.$$

Now suppose that $\zeta_3 \leq \zeta_1 + 1$ (and, hence, $\zeta_2 \in [\zeta_1, \zeta_1 + 1]$). Since $\tilde{\Phi}(\zeta_3 - \zeta_1) \geq \tilde{\Phi}(1)$, it follows that $\zeta_3$ is smaller than the solution to

$$\tilde{\Phi}(1) = 4w\tilde{m}(w)\pi^{-1} \asymp \zeta^{\kappa-1}\frac{g(\zeta)}{\beta(\zeta)}\pi^{-1} \asymp \zeta^{\kappa-1}\phi(\zeta)\pi^{-1}.$$

Taking logarithms, the equation becomes

$$\zeta^2/2 - (\kappa-1)\log\zeta + \log c = \log\pi^{-1},$$

from which it follows that

$$\zeta_2^2 \leq 2\log\pi^{-1} + \log\log\pi^{-1} + C.$$

Consequently, since $\pi[2\log\pi^{-1} + \log\log\pi^{-1} + C]^{q/2}$ is increasing in $\pi$ for sufficiently small $\pi$, and $\pi \leq \eta^{\tilde{p}}\zeta_1^{-p}$, we get

$$\pi\zeta_2^q \leq C\eta^{\tilde{p}}\zeta_1^{-p}[2\log(\eta^{-\tilde{p}}\zeta_1^p)]^{q/2} \leq C\eta^{\tilde{p}}(2\log\eta^{-\tilde{p}})^{q/2}(2\log\tilde{\eta}^{-\tilde{p}})^{-p/2}.$$

If $q < p$, the right-hand side may be bounded further by $C\eta^q$. If $q \geq p$, consider separately the two cases $\eta^{\tilde{p}} > n^{-1}\log^2 n$ and $\eta^{\tilde{p}} \leq n^{-1}\log^2 n$. In all cases we obtain (119) for sufficiently small $\eta$ and $n > n_0$.

To complete the proof of Theorem 1, combine the bounds (112) for $R_q(1)$ and (118) for $\tilde{R}_q(1)$. For the adaptivity bound, similarly combine bounds (111) for $R_q(\zeta_1)$ and (119) for $\tilde{R}_q(\zeta_1)$.

## 9. Proof of Theorem 2.
The proof of Theorem 2 requires small but significant modifications to the proof of Theorem 1.

Consider first the case $\eta^{\tilde{p}} > n^{-1}\log^2 n$, so that $\tilde{\eta} = \eta$. To show that parts (a) and (b) of the theorem remain true with $\hat{\mu}_A$ in place of $\hat{\mu}$, simply observe that

$$E|\delta(X, \hat{t}_A) - \mu|^q = E\{|\delta(X, \hat{t}) - \mu|^q, \hat{t} \leq t_n\} + E\{|\delta(X, t_A) - \mu|^q, \hat{t} > t_n\}.$$

Ignoring the event $\{\hat{t} \leq t_n\}$ in the first term leads to

$$R_q(\hat{\mu}_A, \mu) \leq R_q(\hat{\mu}, \mu) + S_q(\hat{\mu}_A^F, \mu),$$

where

$$(123) \qquad S_q(\hat{\mu}_A^F, \mu) = n^{-1}\sum_i E\{|\delta(X_i, t_A) - \mu_i|^q, \hat{t} > t_n\}.$$

The superscript $F$ emphasizes the fixed threshold $t_A$.

The bound of Theorem 1 applies to $R_q(\hat{\mu}, \mu)$, so it remains to consider $S_q(\hat{\mu}_A^F, \mu)$. Analogously to (100), decompose (123) according to terms with large and small values of $\mu_i$, obtaining

$$(124) \qquad S_q(\hat{\mu}_A^F, \mu) = S_q(\tau) + \tilde{S}_q(\tau),$$

where, for example,

$$\tilde{S}_q(\tau) = n^{-1} \sum_{|\mu_i| \geq \tau} E\{|\delta(X_i, t_A) - \mu_i|^q, \hat{t} > t_n\}.$$

We will need the following risk bounds from Section 6. First, from (71),

$$(125) \qquad E|\delta(X, t_A) - \mu|^q \leq c_q\{|\mu|^q + t_A^{q-1}\phi(t_A)\},$$

while from (70), for any event $B$,

$$(126) \qquad E\{|\delta(X, t_A) - \mu|^q, B\} \leq 4(t_A^q + b^q + 1)P(B)^{1/2}.$$

Consider first part (a), namely, global boundedness. For small $\mu_i$ one uses (125) to obtain

$$(127) \qquad \begin{aligned} S_q(1) &\leq n^{-1} \sum_{|\mu_i| \leq 1} E|\delta(X_i, t_A) - \mu_i|^q \\ &\leq c_q n^{-1} \sum_{|\mu_i| \leq 1} \{|\mu_i|^q + t_A^{q-1}\phi(t_A)\} \\ &\leq c_q\{1 + t_A^{q-1}\phi(t_A)\} \leq 2c_q. \end{aligned}$$

For large $\mu_i$, as in Section 8.4, introduce $\pi = \tilde{\pi}(1, \mu)$ and $\zeta(\mu)$ defined as $\zeta(1, \pi)$, where $\zeta(\tau, \pi)$ is as defined before Lemma 11. Note throughout that $\zeta(\mu) \geq \beta^{-1}(1) > 0$. Arguing as at (117), we also observe that

$$\zeta(\mu) \leq 2\log \pi^{-1}.$$

Two cases arise. If $\zeta(\mu) > \log^{1/2} n$, then

$$t_A = \sqrt{2(1+A)\log n} \leq c\zeta(\mu) \leq c\log \pi^{-1},$$

and so, from (126),

$$E|\eta(X_i, t_A) - \mu_i|^q \leq c(\log \pi^{-1})^q$$

and hence

$$\tilde{S}_q(1) \leq c\pi(\log \pi^{-1})^q \leq C.$$

In the second case $\zeta(\mu) \leq \log^{1/2} n \leq t_n$ and so, using the property $\hat{\zeta} = \zeta(\hat{t}) > \hat{t}$ and Lemma 11,

$$\begin{aligned} P(\hat{t} > t_n) &\leq P(\hat{\zeta} > t_n) \leq P\{\hat{\zeta} > \zeta(\mu)\} \\ &\leq \exp[-Cn\zeta(\mu)^{k-1}\phi\{\zeta(\mu)\}] \leq \exp(-Cn^{1/4}). \end{aligned}$$

Consequently, using (126) with $B = \{\hat{\tau} > t_n\}$,

$$\tilde{S}_q(1) \leq 4\pi(t_A^q + b^q + 1)\exp(-Cn^{1/4}) \leq c\exp(-Cn^{1/4})\log n \leq C.$$

Now turn to part (b), adaptivity over $\ell_p[\eta]$. The case $q \le p$ is simple; from (125) and bound (110), we have

$$S_q(\hat{\mu}_A^F, \mu) \le c_q n^{-1} \sum_i \{|\mu_i|^q + t_A^{q-1} \phi(t_A)\} \le c_q \eta^q + C n^{-(1+A)} \log^{1/2} n.$$

For $q > p$, we follow a strategy broadly similar to that of Section 8.4. In (124) we take $\tau = \zeta_1 = \zeta(\eta, n)$, the pseudothreshold defined by (101). Applying (125) in a similar manner to (127), we find that

$$S_q(\zeta_1) \le c_q n^{-1} \sum_{|\mu_i| \le \zeta_1} \{|\mu_i|^q + t_A^{q-1} \phi(t_A)\} \le c_q \eta^{\widetilde{p}} \zeta_1^{q-p} + c_q t_A^{q-1} \phi(t_A),$$

which is bounded by the right-hand side of (19) in view of (104).

To bound the large signal term $\tilde{S}_q(\zeta_1)$, we again apply Lemma 11 with $\tau = \zeta_2 = \zeta_2(\mu)$ defined as in (120). We first observe, using (126) with $B = \{\hat{t} > t_n\}$, that

$$(128) \qquad \tilde{S}_q(\zeta_1) \le c\pi t_A^q \{P(\hat{t} > t_n)\}^{1/2} \le c\eta^{\widetilde{p}} \zeta^{-p} t_A^q \{P(\hat{t} > t_n)\}^{1/2}.$$

Consider now three cases. First suppose that $\mu$ is such that $\zeta_2(\mu) < t_n$. Using initially the property that $\hat{\zeta} = \zeta(\hat{t}) > \hat{t}$, and then appealing to Lemma 11, we have

$$P(\hat{t} > t_n) \le P(\hat{\zeta} > t_n) \le P(\hat{\zeta} > \zeta_2)$$
$$\le \exp\{-Cn\zeta_2^{\kappa-1} \phi(\zeta_2)\} \le \exp\{-Cnt_n^{\kappa-1} \phi(t_n)\}.$$

Using the definition of $t_n$, and the fact that $t_n^{\kappa-1} \ge 1$ for $n \ge 13$,

$$nt_n^{\kappa-1} \phi(t_n) \ge \phi(0) \log^{5/2} n.$$

Hence, from (128) and using the fact that $\eta \le \eta_0$,

$$\tilde{S}_q(\zeta_1) \le C\eta^{\widetilde{p}} (\log n)^{q/2} \exp(-C \log^{5/2} n) \le C\eta^{\widetilde{p}} \le Cr_{p,q}(\eta).$$

Second, consider $\mu$ for which $\zeta_1^2 \ge \log n$. In this case $2 \log \eta^{-\widetilde{p}} \asymp \zeta_1^2 \ge \log n$, and so, from (128)

$$\tilde{S}_q(\zeta_1) \le C\eta^{\widetilde{p}} (\log n)^{(q-p)/2} \le Cr_{p,q}(\eta).$$

Finally, suppose both $\zeta_1^2 \le \log n$ and $\zeta_2(\mu) \ge t_n$. In this case

$$\zeta_3 - \zeta_1 \ge \zeta_2 - \zeta_1 \ge t_n - \log^{1/2} n \ge \tfrac{1}{2} \log^{1/2} n$$

if $n \ge n_0$. We use (122) defining $\zeta_3$ to derive an upper bound on $\pi$. The equation implies

$$\tilde{\Phi}(\zeta_3 - \zeta_1) = 4\pi^{-1} \tilde{m}(w_3)/\beta(\zeta_3)$$
$$\asymp 4\pi^{-1} \zeta_3^{\kappa-1} g(\zeta_3)/\beta(\zeta_3)$$
$$\asymp 4\pi^{-1} \zeta_3^{\kappa-1} \phi(\zeta_3).$$

In other words, using (103),

$$\pi \asymp \zeta_3^{\kappa-1}(\zeta_3 - \zeta_1)\phi(\zeta_3)/\phi(\zeta_3 - \zeta_1)$$
$$\leq \zeta_3^\kappa \exp\{-(\zeta_3 - \zeta_1)\zeta_1\}\phi(\zeta_1)$$
$$\asymp \eta^p \zeta_1^{1-p} \zeta_3^\kappa \exp\{-(\zeta_3 - \zeta_1)\zeta_1\}.$$

Using the first inequality of (128), we have

$$\tilde{S}_q(\zeta_1) \leq C\eta^{\widetilde{p}} t_A^q \zeta_3^{\kappa+1} \exp\{-(\zeta_3 - \zeta_1)\zeta_1\}$$
$$\leq C\eta^{\widetilde{p}} t_A^q (\zeta_3 - \zeta_1)^{\kappa+1} \exp\{-(\zeta_3 - \zeta_1)\zeta_1\},$$

where we have used the fact that $\zeta_3 - \zeta_1 \geq \frac{1}{2}\log^{1/2} n \geq \frac{1}{2}\zeta_1$, so that $\zeta_3 \leq 3(\zeta_3 - \zeta_1)$. Using these properties again, as well as the property that $\zeta_1 \geq \beta^{-1}(1)$, we have

$$\tilde{S}_q(\zeta_1) \leq C\eta^{\widetilde{p}}(\log n)^{(\kappa+1+q)/2}\exp(-\zeta_1\sqrt{\log n}) \leq C\eta^{\widetilde{p}} \leq Cr_{p,q}(\eta).$$

This completes the proof that the results of Theorem 1 continue to hold for the modified estimator for $\eta^{\widetilde{p}} > n^{-1}\log^2 n$.

Now turn to the case $\eta^{\widetilde{p}} \leq n^{-1}\log^2 n$. Reuse decomposition (100) with $\tau = \zeta(\eta, n)$ defined after (101). First use bound (71) with $t = t_A$:

$$(129) \quad E|\hat{\mu}_{A,i} - \mu_i|^q \leq c_q[|\mu_i|^q + t_A^{q-1}\phi(t_A) + (t_A^q + b^q + 1)\{P(\hat{t}_A < t_A)\}^{1/2}].$$

By the definition of $t_A$, we have

$$(130) \quad t_A^{q-1}\phi(t_A) = \phi(0)n^{-1-A}[2(1+A)\log n]^{(q-1)/2} \leq Cn^{-1-A}\log^{(q-1)/2} n.$$

To bound $P(\hat{t}_A < t_A)$, observe from (53) that $t^2(\zeta) \geq \zeta^2 - C$. In combination with (105), this implies, for $\eta^{\widetilde{p}} \leq n^{-1}\log^2 n$, that

$$t^2(\zeta) \geq t_n^2 + p\log\log n - c - C,$$

so that $t^2(\zeta) \geq t_n^2$ for $n > n(p, \gamma)$. Consequently,

$$\{\hat{t}_A < t_A\} = \{\hat{t} < t_n\} \subset \{\hat{t} < t(\zeta)\} = \{\hat{\zeta} < \zeta\}$$

and so we conclude from (106) that when $\eta^{\widetilde{p}} \leq n^{-1}\log^2 n$ and $n > n(p, \gamma)$,

$$(131) \quad (t_A^q + b^q + 1)P(\hat{t}_A < t_A) \leq c(\log n)^{q/2}\exp\{-C(\log n)^{3/2}\} = o(n^{-1-A}).$$

We now have

$$(132) \quad n^{-1}\sum_i |\mu_i|^q \leq \begin{cases} \left(n^{-1}\sum_i |\mu_i|^p\right)^{q/p} \leq \eta^q = r_{p,q}(\eta), & \text{for } q \leq p, \\ n^{-1}\|\mu\|_p^q \leq n^{(q-p)/p}\eta^q, & \text{for } q > p > 0. \end{cases}$$

Averaging (129) over all $i$ and inserting (130)–(132) proves (23) for the case $q \leq p$ and (24) for $q > p$.

To prove (23) for $q > p$, argue as in Section 8.2 to give

$$n^{-1} \sum_{|\mu_i| \leq \zeta} |\mu_i|^q \leq C r_{p,q}(\eta), \tag{133}$$

so that, summing only over $|\mu_i| \leq \zeta$,

$$R_q(\zeta) \leq C r_{p,q}(\eta) + c n^{-1-A} \log^{1/2} n. \tag{134}$$

For $q > p$ and $|\mu_i| > \zeta$, we apply bound (78), noting that $\hat{t}_A \leq t_A$ with probability one, to obtain

$$E|\hat{\mu}_{A,i} - \mu_i|^q \leq 4(t_A^q + b^q + 1) \leq C(\log n)^{q/2}.$$

As in the previous section, comparing (120) and (121), we have

$$\pi = n^{-1} \#\{i : |\mu_i| > \zeta\} \leq \eta^{\widetilde{p}} \zeta^{-p},$$

and so, recalling from (105) that $\zeta > \sqrt{\log n}$,

$$\tilde{R}_q(\zeta) \leq \pi C (\log n)^{q/2} \leq C \eta^{\widetilde{p}} (\log n)^{(q-p)/2}.$$

But $\eta^{\widetilde{p}} \leq n^{-1} \log^2 n$ implies that $\log \eta^{-\widetilde{p}} \geq \log n - 2 \log \log n$ and, hence, that $\log n \leq C \log \eta^{-\widetilde{p}}$, so for $n$ large and $\eta^{\widetilde{p}} \leq n^{-1} \log^2 n$ we have

$$\tilde{R}_q(\zeta) \leq C r_{p,q}(\eta);$$

combining this result with (134) completes the proof of Theorem 2.   $\square$

**10. Remarks on the posterior mean.** In proving results for the posterior mean, we have assumed throughout that $q > 1$. The failure of the posterior mean to be a strict thresholding rule has a substantive effect on the overall risk if $q \leq 1$. Concentrate attention on the case where $\gamma$ is the Laplace distribution with parameter 1, and define $\psi(n) = \exp(\sqrt{2 \log n})$ so that $g(\sqrt{2 \log n})^{-1} \asymp \psi(n)$ as $n \to \infty$.

An important contributor to our arguments was $t(\hat{w}) \leq \sqrt{2 \log n}$, from which it follows that $\tau(\hat{w}) \leq \sqrt{2 \log n}$. By the definition of $\tau$, we have

$$\frac{\hat{w}}{1 - \hat{w}} \geq \frac{\phi(\sqrt{2 \log n})}{g(\sqrt{2 \log n})} \geq C n^{-1} \psi(n),$$

so that $\hat{w} \geq C n^{-1} \psi(n)$. Since $g(x)/\phi(x)$ is bounded below away from zero, it follows that, for some constant $C$ and for all $x$, the posterior weight $\tilde{w}(x, \hat{w}) \geq C n^{-1} \psi(n)$.

On the other hand, the odd function $\tilde{\mu}_1$ satisfies $\tilde{\mu}_1'(0) > 0$ and is strictly increasing with $\tilde{\mu}_1(x) \geq x - \Lambda$ for all $x$; therefore, $x^{-1}\tilde{\mu}_1(x)$ is uniformly bounded below away from zero for $x \neq 0$. It follows that, for all $x \neq 0$,

$$|\tilde{\mu}(x, \hat{w})| = \tilde{w}(x, \hat{w})|\tilde{\mu}_1(x)| \geq C|x|n^{-1}\psi(n).$$

If $\mu = 0$ and $X \sim N(0, 1)$, it follows that

$$E|\tilde{\mu}(X, \hat{w}) - \mu|^q \geq C^q n^{-q}\psi(n)^q E|X|^q = Cn^{-q}\psi(n)^q$$

so that, however small the value of $\eta$, the risk bound cannot be reduced below $Cn^{-q}\psi(n)^q$, making it impossible for the estimate to attain the full range of adaptivity given by the posterior median. The restrictions become more severe the lower the value of $q$.

**Acknowledgments.** I. M. Johnstone is very grateful for the hospitality of the University of Bristol, the "Nonparametric Semester" of Institut Henri Poincaré, Paris, and the Australian National University, where parts of the work on this paper were carried out. B. W. Silverman acknowledges the hospitality of the Department of Statistics and the Center for Advanced Study in the Behavioral Sciences at Stanford University. An early version of this work formed part of his Special Invited Paper presented in 1999. We thank Tim Heaton for carefully reading the manuscript and pointing out a number of minor errors. We owe a particular debt of gratitude to Ed George for his encouragement and generosity over the whole period of this work. The careful suggestions of the referees were extremely useful in prompting us to improve both the presentation and the results of the paper.

DEPARTMENT OF STATISTICS
STANFORD UNIVERSITY
STANFORD, CALIFORNIA 94305-4065
USA
E-MAIL: imj@stat.stanford.edu

ST. PETER'S COLLEGE
OXFORD OX1 2DL
UNITED KINGDOM
E-MAIL: silverman@stats.ox.ac.uk